\pdfoutput=1
\documentclass{article}
\usepackage{hyperref} 
\usepackage{amsthm}
\usepackage{amssymb}
\usepackage{graphicx}
\usepackage{color}
\usepackage{framed}   
\usepackage{color}
\usepackage{tikz}
\usetikzlibrary{calc}
\usepackage{amstext}
\usepackage{amssymb}
\usepackage{amsmath}
\usepackage{graphicx}
\usepackage{graphics} 
\usepackage[all,arc,poly]{xy}
\usepackage{color}
\usepackage{relsize}

\renewcommand{\phi}{\varphi}

\newcommand{\nc}{\newcommand}
\nc{\Section}{\section}
\nc{\SubSection}{\subsection}
\nc{\grey}[1]{\textcolor{lightgray}{#1}}

\newtheorem{theo}{Theorem}[section]   
\newtheorem{ddef}[theo]{Definition}
\newtheorem{llem}[theo]{Lemma} 
\newtheorem{oobs}[theo]{Observation} 
\newtheorem{rrem}[theo]{Remark} 
\newtheorem{prop}[theo]{Proposition} 
\newtheorem{ccor}[theo]{Corollary}
\newtheorem{claim}[theo]{Claim}  
\newtheorem{qquest}[theo]{Question} 
\newtheorem{fact}[theo]{Fact} 
\newtheorem{pprov}[theo]{Proviso}
\newtheorem{eexam}[theo]{Example} 

\nc{\bT}{\begin{theo}} 
\nc{\eT}{\end{theo}}
\nc{\bD}{\begin{ddef} \rm }
\nc{\eD}{\end{ddef}}
\nc{\bC}{\begin{ccor}}
\nc{\eC}{\end{ccor}}
\nc{\bCl}{\begin{claim}}
\nc{\eCl}{\end{claim}}
\nc{\bQ}{\begin{qquest}}
\nc{\eQ}{\end{qquest}}
\nc{\bL}{\begin{llem}}
\nc{\eL}{\end{llem}}
\nc{\bP}{\begin{prop}}
\nc{\eP}{\end{prop}}
\nc{\bR}{\begin{rrem}}
\nc{\eR}{\end{rrem}}
\nc{\bO}{\begin{oobs}}
\nc{\eO}{\end{oobs}}
\nc{\bF}{\begin{fact}}
\nc{\eF}{\end{fact}}
\nc{\bProv}{\begin{pprov}}
\nc{\eProv}{\end{pprov}}
\nc{\bE}{\begin{eexam} \rm }
\nc{\eE}{\end{eexam}}

\renewcommand{\geq}{\geqslant}
\renewcommand{\leq}{\leqslant}

\renewcommand{\subset}{\subseteq}

\newcommand{\kbox}{\ensuremath{\square}}
\newcommand{\ibox}{\boxplus}
\def\disji{\rotatebox[origin=c]{-90}{$\!{\geqslant}$}}
\newcommand{\lori}{\,\disji\,}
\newcommand{\simn}{{\sim^n}}
\newcommand{\simm}{{\sim^m}}

\newcommand{\ee}{\varepsilon}

\nc{\simg}{\sim_{\mathsf{g}}}
\nc{\ssimg}{\approx_{\mathsf{g}}}
\nc{\FO}{\mathsf{FO}}
\nc{\MSO}{\mathsf{MSO}}
\nc{\ML}{\mathsf{ML}}
\nc{\IML}{\mathsf{IML}}
\nc{\SO}{\mathsf{SO}}
\nc{\GF}{\mathsf{GF}}
\nc{\gd}{\mathsf{gd}}
\nc{\hg}{\mathsf{hg}}
\nc{\vg}{\mathsf{vg}}
\nc{\free}{\mathrm{free}}
\nc{\dom}{\mathrm{dom}}
\nc{\nn}{\mathsf{non}}
\nc{\NE}{\mathsf{NE}}
\nc{\NF}{\mathsf{NF}}
\nc{\C}{\mathrm{cl}}
\nc{\expl}[1]{\mathsf{I}(#1)}
\nc{\INV}[1]{\mathsf{Inv}(#1)} 

\nc{\full}{\mathsf{full}}
\nc{\lf}{\mathsf{lf}}
\nc{\dc}{\mathsf{dc}}
\nc{\rel}{\mathsf{rel}}

\newenvironment{romanenumerate}%
{\begin{list}{(\roman{enumi})}{\usecounter{enumi}
\setlength{\labelwidth}{2cm}
\setlength{\itemindent}{0pt}
\setlength{\itemsep}{0.5\itemsep}
\setlength{\topsep}{\itemsep}
\setlength{\parsep}{0pt}
}}{\end{list}}
\nc{\bre}{\begin{romanenumerate}}
\nc{\ere}{\end{romanenumerate}}
\newenvironment{alphaenumerate}%
{\begin{list}{(\alph{enumii})}{\usecounter{enumii}
\setlength{\labelwidth}{2cm}
\setlength{\itemindent}{0pt}
\setlength{\itemsep}{0.5\itemsep}
\setlength{\topsep}{\itemsep}
\setlength{\parsep}{0pt}
}}{\end{list}}
\nc{\bae}{\begin{alphaenumerate}}
\nc{\eae}{\end{alphaenumerate}}
\newenvironment{numenumerate}%
{\begin{list}{(\arabic{enumiii})}{\usecounter{enumiii}
\setlength{\labelwidth}{2cm}
\setlength{\itemindent}{0pt}
\setlength{\itemsep}{0.5\itemsep}
\setlength{\topsep}{\itemsep}
\setlength{\parsep}{0pt}
}}{\end{list}}
\nc{\bne}{\begin{numenumerate}}
\nc{\ene}{\end{numenumerate}}

\nc{\ins}[1]{\bigskip\noindent
\framebox{\begin{minipage}{.98\textwidth} \sloppy \noindent \em #1 \end{minipage}}\bigskip}

\nc{\str}[1]{{\mathfrak{#1}}}
\nc{\brck}[1]{[\![ #1 ]\!]}
\nc{\restr}{\!\restriction\!}

\nc{\HH}{\mathbb{H}}
\nc{\VV}{\mathbb{V}}
\nc{\abar}{\mathbf{a}}
\nc{\bbar}{\mathbf{b}}
\nc{\cbar}{\mathbf{c}}
\nc{\mbar}{\mathbf{m}}
\nc{\xbar}{\mathbf{x}}
\nc{\ybar}{\mathbf{y}}
\nc{\zbar}{\mathbf{z}}
\nc{\ubar}{\mathbf{u}}
\nc{\sbar}{\mathbf{s}}
\nc{\tbar}{\mathbf{t}}
\nc{\vbar}{\mathbf{v}}
\nc{\wbar}{\mathbf{w}}

\nc{\Xbar}{\mathbf{X}}
\nc{\Ybar}{\mathbf{Y}}
\nc{\Zbar}{\mathbf{Z}}
\nc{\Pbar}{\mathbf{P}}
\nc{\nubar}{\mbox{\boldmath $\nu$}}

\nc{\barr}{\begin{array}}
\nc{\earr}{\end{array}}
\nc{\btab}{\begin{tabular}}
\nc{\etab}{\end{tabular}}

\nc{\nothing}{\rule{0em}{1ex}}
\nc{\highnothing}{\rule{0em}{3ex}}
\nc{\hnt}{\highnothing}
\nc{\nt}{\nothing}
\nc{\nnt}{\rule{.1pt}{0pt}}

\nc{\scc}{\scriptstyle}
\nc{\ssc}{\scriptscriptstyle}
\nc{\N}{{\mathbb N}}
\nc{\Z}{{\mathbb Z}}
\nc{\M}{{\mathbb M}}
\nc{\F}{{\mathbb F}}
\nc{\W}{{\mathbb W}}
\newcommand{\PI}{\mbox{\bf I}}
\newcommand{\PII}{\mbox{\bf II}}
\renewcommand{\P}{\ensuremath{\mathcal{P}}}
\newcommand{\A}{\ensuremath{\mathcal{A}}}
\newcommand{\CC}{\ensuremath{\mathcal{C}}}

\newcommand{\MM}{\ensuremath{\mathfrak{M}}}

\newcommand{\w}{\texttt{w}}
\renewcommand{\v}{\texttt{v}}

\newcommand{\inqbm}{\textsc{InqML}}
\newcommand{\inqml}{\textsc{InqML}}

\nc{\prf}{\begin{proof}}
\nc{\eprf}{\end{proof}}
 
\usepackage{wrapfig} 

\tikzstyle{index on}=[inner sep=2pt, white, circle, fill=black]
\tikzstyle{index off}=[inner sep=2pt, black, circle, draw]
\tikzstyle{index gray}=[inner sep=2pt, black, circle, fill=lightgray]
\tikzstyle{opaque}=[fill=gray,fill opacity=.1]

\title{On the expressive power of \\ inquisitive epistemic logic} 

\author{Ivano Ciardelli and Martin Otto}
\date{}

\begin{document}
\maketitle

\begin{abstract}
  \noindent
Inquisitive modal logic, \inqml, in its epistemic incarnation, extends
standard epistemic logic to capture not just the information that
agents have, but also the questions that they are interested in.
We use the natural notion of bisimulation equivalence  in the setting
of \inqml, as introduced in~\cite{CiardelliOttoJSL2021}, to
characterise the expressiveness of \inqml\ as the bisimulation
invariant fragment of first-order logic over 
natural classes of two-sorted first-order structures that
arise as relational encodings of inquisitive epistemic
(S5-like) models.
The non-elementary nature of these classes
crucially requires non-classical model-theoretic methods for the
analysis of first-order expressiveness,
irrespective of whether we aim for characterisations in the sense of classical 
or of finite model theory. 
\end{abstract}

\pagebreak

\tableofcontents

\pagebreak

\section{Introduction}
\label{introsec}

Inquisitive logic \cite{CiardelliRoelofsen,Ciardelli:23book} is a framework that allows us to extend various systems of classical logic with formulae that represent questions. Thus, in the propositional setting, in addition to the formula $p$ representing the statement that $p$, we also have, e.g., a formula $?p$ representing the question \emph{whether $p$}. In the first-order setting, in addition to the formula $\forall xPx$ representing the statement that all objects are $P$, we also have, e.g., a formula $\forall x?Px$, representing the question \emph{which objects are $P$}. Technically, this extension is achieved by using standard models for intensional logic, but interpreting formulae not relative to single possible worlds, but relative to sets of worlds, called \emph{information states}. As such, inquisitive logics fit within the family of logics based on versions of \emph{team semantics}, where formulae are evaluated relative to sets of (standard) semantic parameters.\footnote{Team semantics has been used to interpret independence-friendly logic \cite{Hodges:97,Hodges:97b} and many variants of dependence and independence logic (see, among others, \cite{Vaananen:07,Vaananen:08,AbramskyVaananen:09,Galliani:12,Yang:14,YangVaananen}).}

Inquisitive epistemic logic
\cite{CiardelliRoelofsen:15idel,Ciardelli:14aiml,IvanoDiss,Gessel:20action,PuncocharSedlar:21epistemic}
is an extension of standard epistemic logic~\cite{Hintikka:62}, 
which, building on inquisitive logic, allows us to reason not only about the knowledge that agents have, but also about the issues they are interested in. In this logic we can, for instance, embed the question $?p$ under suitable inquisitive modalities, resulting in modal statements which express that a certain agent \emph{knows whether $p$}, or that a certain agent \emph{wonders whether~$p$}. Furthermore, in a dynamic extension of the logic, we can model not only the communicative effect of publicly making statements (as in public announcement logic \cite{Plaza:07}), but also the effect of publicly asking questions, describing in particular how both kinds of moves affect the issues agents are interested in.

From a technical point of view, inquisitive epistemic logic is a special case of inquisitive (multi-)modal logic. As in Kripke semantics, models for this logic assign to each world $w$ and agent $a$ a set $\Sigma_a(w)$ of successors, which inquisitive modalities then quantify over; however, differently from Kripke semantics, these successors are not worlds, but rather information states, i.e., sets of worlds. This means that, in the inquisitive setting, modalities involve a restricted kind of second-order quantification. Within this general modal setup, inquisitive epistemic logic arises as a special case when we impose $S5$-like constraints on the modal frame, motivated as usual by the assumptions of factivity of knowledge and full introspection of agents.

Our topic in this paper is the expressive power of inquisitive epistemic logic, viewed as a language for describing properties of worlds as well as properties of information states in the setting of inquisitive epistemic models. In particular, we are interested in comparing the expressive power of inquisitive epistemic logic with that of classical first-order predicate logic $\FO$ interpreted over (suitable relational encodings of) inquisitive epistemic models. 

A crucial role in this study is played by a notion of \emph{bisimilarity} for inquisitive modal models, identified in our previous paper~\cite{CiardelliOttoJSL2021}. For inquisitive modal logic interpreted over arbitrary (or arbitrary finite) models, a van Benthem/Rosen characterization result was proved in \cite{CiardelliOttoJSL2021}.

\bT[\cite{CiardelliOttoJSL2021}]
\label{main1}
Inquisitive modal logic can be characterised as the 
bi\-si\-mu\-la\-tion-invariant fragment of first-order logic $\FO$
over some natural classes of (finite or arbitrary) relational inquisitive models. 
\eT

In this paper, we will prove an analogous result for the specific setting of inquisitive epistemic logic, where models are required to obey $S5$-like constrains.

\bT[main theorem]
\label{main2}
Inquisitive epistemic logic (in a multi-agent setting) can be characterised as the 
bisimulation-invariant fragment of $\FO$ over some natural
classes of (finite or arbitrary) relational inquisitive epistemic models. 
\eT

Extending Theorem~\ref{main1} to the inquisitive epistemic case is a challenging enterprise:
the reason is that the proof in \cite{CiardelliOttoJSL2021} crucially relies on a 
step in which we unfold a model into a bisimilar locally tree-like model. 
That step, however, is incompatible with the $S5$-like requirements on
inquisitive epistemic models, which rule out simple tree-like structures. 
We here show that locally tree-like behaviour
can be achieved for relational structures in which suitable bisimilar
companions are first-order interpretable. Our analysis of the
bisimulation types of epistemic inquisitive models
and their relational encodings
is then based on an analysis 
of these representations that supports the central upgrading argument 
linking finite levels of bisimulation equivalence to levels of
elementary equivalence between the interpreted models.
Equivalences between  the locally tree-like representations are
assessed through a combination of classical back-and-forth techniques
that capture the expressive power of first-order as well as of restricted
forms of monadic second-order logic.
In these aspects the current paper technically builds on but also considerably exceeds
the treatment in~\cite{CiardelliOttoJSL2021}.

Besides its technical interest, 
Theorem \ref{main2} is a significant result from a conceptual point of view. It entails that inquisitive epistemic logic is a natural choice for a language designed to express bisimulation-invariant properties of pointed inquisitive epistemic models. On the one hand, this language allows us to express \emph{only} properties that are bisimulation invariant;
this is desirable if one views bisimilar points as distinct
representations of what is essentially the same state of affairs, so
that aspects of a model which are not bisimulation invariant are spurious
artefacts of the representation that our logical language should be insensitive to. On the other hand, as the theorem states, this language allows us to express \emph{all} the bisimulation-invariant properties that are expressible in first-order predicate logic, thus enjoying a form of expressive completeness.

\medskip
The paper is structured as follows: Section~2 introduces all the
relevant background for our contribution; Section~3 contains the steps
leading up to the proof of Theorem \ref{main2};
Section~4 summarises our results and outlines some directions for further work.

\newpage
\section{Background}
\label{sec:inquisitive modal logic}

In this section we review some background for our contribution. First,
we briefly introduce inquisitive epistemic logic
\cite{CiardelliRoelofsen:15idel,IvanoDiss}, as well as the associated
notion of bisimulation \cite{CiardelliOttoJSL2021}. We also show how
to encode models of inquisitive epistemic models as relational
structures, following \cite{CiardelliOttoJSL2021},  and recall how the
expressive
power of first-order logic and monadic second-order logic can be characterised by means of Ehrenfeucht--Fra\"\i ss\'e games.

\subsection{Inquisitive epistemic models}

Let us fix a \emph{signature} for our logic, consisting of a set \A\ of agents and a set \P\ of atomic propositions. Relative to such a signature, the notions of inquisitive epistemic frames and models are defined as follows. 

\bD[inquisitive epistemic models]
\label{epinqdef}\\
An \emph{inquisitive epistemic frame} is a pair
$\F=( W,(\Sigma_a)_{a\in\A})$, where $W$ is a non-empty set whose elements are called \emph{worlds}, and for each $a\in\A$, the map $\Sigma_a:W\to\wp\wp(W)$ 
assigns to each world
$w$ a set of sets of worlds $\Sigma_a(w)$ in accordance with the 
following constraints, where $\sigma_a(w):=\bigcup\Sigma_a(w)$. 
For all $w,v\in W$, $s,t\subseteq W$, and $a \in \A$ we must have that:
\begin{itemize}
\item[--] $\Sigma_a(w)$ is non-empty;
\item[--] if $s\in \Sigma_a(w)$ and $t\subseteq s$ then $t\in\Sigma_a(w)$ (downward closure);
\item[--] $w\in\sigma_a(w)$ (factivity);
\item[--] if $v\in\sigma_a(w) $ then $\Sigma_a(v)=\Sigma_a(w)$
(full introspection).
\end{itemize}
\noindent 
An inquisitive epistemic model
consists of an inquisitive epistemic frame together with a
propositional assignment $V \colon \P \to\wp(W)$.
\eD

A set of worlds $s\subseteq W$ is referred to
as an \emph{information state}; a non-empty set of information states which is
downward closed (i.e., such that whenever it contains an information state $s$, it contains every subset of $s$) is called an \emph{inquisitive state}. Thus, the above definition says that each map $\Sigma_a$ assigns to a world $w$ an inquisitive state $\Sigma_a(w)$,
in accordance with factivity and introspection.

In the intended interpretation for inquisitive epistemic models \cite{CiardelliRoelofsen:15idel,Ciardelli:14aiml},
the map $\Sigma_a$ is taken to describe both the \emph{knowledge} and the \emph{issues}
of agent~$a$. 
The agent's knowledge state at $w$,
$\sigma_a(w)=\bigcup\Sigma_a(w)$, consists of all the worlds that are
compatible with what the agent knows. The agent's inquisitive state at
$w$, $\Sigma_a(w)$, consists of all those information states where the
agent's issues are settled. This interpretation motivates the first two 
constaints on $\Sigma_a(w)$: $\Sigma_a(w)$ should be downward closed, 
since if the agent's issues are settled by an information state $s$, they
will also be settled by any stronger information state $t\subset s$; 
furthermore, the agent's issues will always be trivially settled by the
empty information state, $\emptyset$, which represents the state
of inconsistent information. The factivity constraint amounts, as usual,
to the requirement that the agent's knowledge at a world be truthful.
The introspection constraint requires agents to have full knowledge
of their own inquisitive state: at each world $w$, only worlds $v$
where the agent's state is the same as in $w$ will be compatible 
with the agent's knowledge. 

Note that with an inquisitive epistemic frame $\F$ we can always associate a
standard Kripke frame $\str{K}(\F)$ having the same set of worlds and
accessibility relations defined by $R_a \!:=\! \{ ( v,w) \colon \!v\! \in\! \sigma_a(w) \}$. 
Similarly, with an inquisitive epistemic model $\M$ we can associate a
standard Kripke model $\str{K}(\M)$.  
It is easy to verify that the Kripke frame associated with an inquisitive
epistemic frame is always an $S5$ frame, i.e., all the accessibility relations 
$R_a$ are equivalence relations on~$W$.
We denote as $[w]_a$ the equivalence class of world $w$ w.r.t.\ $R_a$,
and by $W/R_a = \{ [w]_a\colon w \in W \}$ the quotient consisting of
these equivalence classes, which form a partition of the set $W$ of
possible worlds, reflecting the possible knowledge states of agent~$a$.
The introspection condition ensures that the map $\Sigma_a$ takes
a constant value on each of these equivalence classes.

\subsection{Inquisitive modal logic}

The syntax  of inquisitive modal logic $\inqbm$ is given by:
\[
  \phi \; ::= \;
  p\;|\,\bot\,|\,(\phi\land\phi)\,|\,(\phi\to\phi)\,|\,(\phi\lori\phi)\,|\, \kbox_a\phi\,|\,\ibox_a\!\phi,
\]
where $p \in \P$ stands for a basic proposition, and
$a \in \A$ is the labelling of inquisitive modalities by agents. 
Nesting depth w.r.t.\ the modal operators $\kbox_a$ and $\ibox_a$
(combined and across all $a \in \A$) is syntactically defined as usual, and for $n
\in \N$ we let
$\inqbm_n \subset \inqbm$ stand for the restriction to nesting depth
up to~$n$.
The mono-modal case, which is the focus of \cite{CiardelliOttoJSL2021}, 
may be regarded as a special case, where the label for a single agent is
omitted as irrelevant.

We treat negation and disjunction as defined connectives (syntactic shorthands) 
according to:
$$\neg\phi:=\phi\to\bot\qquad\qquad\phi\lor\psi:=\neg(\neg\phi\land\neg\psi)$$
With these abbreviations, the above syntax can be seen as including
the syntax of standard
modal logic, based on Boolean connectives and the modalities $\Box_a$. 
In particular, $\inqbm_0$ includes all formulae standard propositional logic. 

In addition to standard connectives, the language contains
\emph{inquisitive disjunction} $\lori$, which is interpreted as
a question-forming (inclusive) disjunction: for instance, whereas
the formula $p\lor\neg p$ stands for the (tautological) statement 
that $p$ is either true or false, the formula $p\lori\neg p$ (abbreviated 
as $?p$) stands for the question \emph{whether} $p$ is true or false.
Equally importantly, the language contains two modalities for each
agent, which allow us to embed both statements and questions. Both
kinds, $\kbox_a$ and $\ibox_a$ coincide with a standard Kripke
box modality $\kbox_a$ when applied to statements, but 
differ in application to questions.%
\footnote{Modalities
$\kbox_a$ and $\ibox_a$ are usually read as ``box'' and ``window''. In the epistemic setting (e.g., in ~\cite{CiardelliRoelofsen:15idel}), they are denoted $K_a$ and $E_a$, and
read as ``know'' and ``entertain'' respectively.}

In line with the inquisitive approach, the semantics 
for \inqbm\ is not given by a recursive definition of \emph{truth} relative
to a possible world, but rather by a recursive definition of a notion of
\emph{support} relative to an information state, i.e., a set of
worlds.

\bD[semantics of $\inqbm$]~\\
\label{supportsemdefn}%
Let $\M=(W,(\Sigma_a),V)$ be an inquisitive epistemic model, $s\subseteq W$:
\begin{itemize}
\item $\M,s\models p\iff s\subseteq V(p)$
\item $\M,s\models \bot\iff s=\emptyset$
\item $\M,s\models \phi\land\psi\iff \M,s\models\phi\text{ and }\M,s\models\psi$
\item 
$\M,s\models \phi\to\psi\iff \mbox{ for all } t\subseteq s:\M,t\models\phi\;\Rightarrow\; \M,t\models\psi$
\item $\M,s\models \phi\lori\psi\iff \M,s\models\phi\text{ or }\M,s\models\psi$
\item $\M,s\models \kbox_a\phi\iff \mbox{ for all } w\in s: \M,\sigma_a(w)\models\phi$
\item $\M,s\models \ibox_a\phi\iff \mbox{ for all } w\in s, \mbox{ for all } t\in\Sigma_a(w): \M,t\models\phi$
\end{itemize}
\eD

The following two properties hold generally in \inqbm. 

 \bP\label{prop:persistency}~
 \begin{itemize}
 \item persistency: if $\M,s\models\phi$ and $t\subseteq s$, then $\M,t\models\phi$;
   \item empty state property (semantic ex-falso): $\M,\emptyset\models\phi$ for all $\phi\in\inqbm$.
 \end{itemize}
 \eP

 The first item says that support is preserved as information
 increases, i.e., as we move from a state to an extension of it. The
 second item says that the inconsistent 
 information state vacuously supports every formula; 
 this can be seen as a semantic counterpart of the 
 \emph{ex falso quodlibet} principle.  
These
 principles imply that the support set $[\phi]_\M := \{ s \subseteq W
 \colon \M,s \models \phi \}$ of a formula is downward closed and
 non-empty, i.e., it is an inquisitive state.

 Although the primary notion in \inqml\ is support at an information state, truth at a world is retrieved as a derivative notion.

 \bD[truth]~\\ 
 \label{truthdef}
 $\phi$ is \emph{true} at a world $w$ in a model $\M$, 
 denoted $\M,w\models\phi$, in case $\M,\{w\}\models\phi$.
 \eD

 Spelling out Definition~\ref{supportsemdefn} in the special case of
 singleton states, we see that standard connectives have the usual
 truth-conditional behaviour. 
 For modal formulae, we find the following truth-conditions:
 \begin{itemize}
 \item $\M,w\models \kbox_a\phi\iff \M,\sigma_a(w)\models\phi$
 \item $\M,w\models \ibox_a\phi\iff \forall t\in\Sigma_a(w): \M,t\models\phi$
 \end{itemize}

 Notice that truth in \inqbm\ cannot be given a direct recursive
 definition, as the truth conditions for modal formulae $\kbox_a\phi$ and
 $\ibox_a\phi$ depend on the support conditions for $\phi$ --- not just on
 its truth conditions.

 For many formulae, support at a state just boils down to truth at each
 world in the state. We refer to these formulae as \emph{truth-conditional}.%
 \footnote{In team semantic terminology 
 (e.g., \cite{Vaananen:07,YangVaananen}),
 truth-conditional formulae are called \emph{flat}.}

 \medskip
  \bD[truth-conditional formulae]\label{flatdef}~\\
  We say that a formula $\phi$ is \emph{truth-conditional} if for all
  models $\M$ and information states $s$: $\M,s\models\phi\iff \M,w\models\phi$ for all $w\in s$.
  \eD

 Following \cite{IvanoDiss}, we view truth-conditional 
 formulae as statements, and non-truth-conditional formulae as
 questions. 
 It is immediate to verify that atomic formulae, $\bot$, 
 and all formulae of the form $\kbox_a\phi$ and $\ibox_a\phi$ are truth-conditional.
 Furthermore, the class of truth-conditional formulae is closed under all
 connectives except for $\lori$. This means that if all occurrences of $\lori$ in $\phi$
 are within the scope of a modality, then $\phi$ is truth-conditional.
 Using this fact, it is easy to see that all formulae of standard modal logic, i.e., formulae which do not contain occurrences of $\lori$ or $\ibox_a$, receive exactly the same truth-conditions as in standard modal logic.

 \bP 
 If $\phi$ is a formula not containing $\lori$ or $\ibox_a$, then $\M,w\models\phi\iff\str{K}(\M),w\models\phi$ in standard Kripke semantics.
 \eP

 As long as questions are not around, the modality~$\ibox_a$ also
 coincides with $\kbox_a$, and with the standard box modality. 
 That is, if $\phi$ is truth-conditional, then 
 \[
 M,w\models\kbox_a\phi\iff M,w\models\ibox_a\phi\iff M,v\models\phi\text{
   for all }v\in\sigma_a(w).
 \]

By contrast, the two modalities come
apart when they are applied to questions. Consider, for instance,
the formulae $\kbox_a{?p}$ and $\ibox_a{?p}$. According to the above support definition, the formula $?p$ (i.e., $p\lori\neg p$)
is supported by an information state $s$ just in case $p$ has the same
truth value at all worlds in $s$---in other words, just in case the information
available in $s$ settles the question whether~$p$.
Thus, according to the truth conditions for modal formulae, 
 $\kbox_a{?p}$ is true at a world $w$ iff the knowledge state $\sigma_a(w)= [w]_a$
of agent~$a$, settles the question $?p$. Thus, $\kbox_a{?p}$ expresses
the fact that the agent \emph{knows whether $p$}. By contrast, $\ibox_a{?p}$ is
true iff any information state $t\in\Sigma_a(w)$ (i.e., any state that
settles the agent's issues) also settles $?p$; thus, roughly speaking, $\ibox_a{?p}$
expresses that settling the question whether $p$ is part of the agent's goals, i.e., the agent \emph{is 
interested in whether $p$}  
(see \cite{CiardelliRoelofsen:15idel,IvanoDiss} for more detailed discussion).

Inquisitive modal formulae can be interpreted both relative to information states and (derivatively) relative to worlds. They can thus be seen both as expressing properties of state-pointed models, and as expressing properties of world-pointed models. In general, by a property of world-pointed or state-pointed models we mean a class of such pointed models. A property is \emph{definable} in $\inqbm$ if it is the class of models 
where some formula $\phi \in \inqbm$ is true (in the
world-pointed case) or supported (in the state-pointed case).%
\footnote{In connection with \inqbm-definability by single
  formulae we implicitly restrict attention to finite modal signatures, i.e.\
\CC\ here is, w.l.o.g.\ specified in a signature with finite $\A$ and $\P$.}
For the case of state-pointed models, Proposition \ref{prop:persistency} implies that any definable property will be an \emph{inquisitive property}, in the sense of the following definition.

\bD[inquisitive properties]\label{inqclass}
A property $\mathcal{K}$ of state-pointed models is said to be \emph{inquisitive}
if it is downward closed with respect to the information state 
(i.e., $(\M,s) \in \mathcal{K}$ implies $(\M,t) \in \mathcal{K}$ for all $t\subseteq s$) and satisfies the empty state property 
 (i.e., $(\M,\emptyset) \in \mathcal{K}$
  for all $\M$).
\eD

In the rest of the paper, when dealing with properties of state-pointed models, we can thus restrict our attention to inquisitive properties.
We now turn to the task of giving a precise characterization of which world-properties and inquisitive state-properties are definable in \inqml.

\subsection{Inquisitive bisimulation}
\label{sec:bisimulations} 

An inquisitive modal model can be seen as a structure with two
closely linked sorts of entities: worlds and information states.
Information states $s$ are sets of worlds, determined
by their elements; worlds $w$ carry, besides the propositional assignment,
inquisitive assignments $\Sigma_a(w)$, which are sets of of information states.
The natural notion of bisimilarity, and its 
finite approximations of $n$-bisimilarity for $n \in \N$, capture
notions of equivalence based on a step-wise probing of these links. 
We here focus on a corresponding game-based formulation for
levels of bisimulation equivalence in terms of winning
strategies in two-player bisimulation games, which highlight the
dynamic nature of the``probing'' of behavioural equivalence.
(The alternative
formalisation in terms of back-and-forth systems is easily seen to be
equivalent as usual, cf.~\cite{CiardelliOttoJSL2021} for more
detail.) 

The game is played by two players, \PI\ and \PII, 
who act as challenger and defender of a similarity claim involving a
pair of worlds $w$ and $w'$ or information states $s$ and $s'$ over two models
$\M=( W,(\Sigma_a),V )$ and $\M'=( W',(\Sigma_a'),V')$. 
We distinguish world-positions $( w,w' ) \in W \times W'$ and 
state-positions as $( s,s' ) \in \wp(W) \times \wp(W')$,
which match two worlds or two information states, respectively.
The game alternates between world-positions and
state-positions, reflecting the links between information states and 
worlds in $\M$ and $\M'$, with a view to testing or probing
similarities or distinguishing features for the currently matched
pairs, as follows. In a world-position $( w,w' )$,
  \PI\ chooses an agent $a \in \A$ and an information state in the
  inquisitive state associated to one of these worlds, i.e.\ either
  some $s \in \Sigma_a(w)$ or some 
$s' \in \Sigma_a'(w')$, and \PII\ must respond by choosing an
information state for the same agent~$a$ on the opposite side; this
exchange results in a state-position $( s,s' )$.  
In a state-position
$( s,s')$, \PI\ chooses a world in either state, i.e.\ either some $w
\in s$ or some $w' \in s'$,  and \PII\ must respond likewise
on the opposite side; this results in a world-position $( w,w' )$.

A round of the game consists of four moves leading from a world-position to another. 
In the bounded version of the game, the number of rounds is fixed
in advance, while the game is allowed to go on indefinitely
in the unbounded version. Either player loses when stuck for a move,
which in epistemic models can only happen in state positions involving
the empty information state. The game also ends, 
with a loss for \PII, in any world-position $( w,w')$
in which $w$ and $w'$ disagree on the truth of some $p\in\P$. 
All other plays, and in particular infinite runs of the unbounded game, 
are won by $\PII$. The following definition applies to arbitrary inquisitive 
modal models (which do not necessarily satisfy factivity and introspection), 
but here we are specifically interested in its application to the case of 
inquisitive epistemic models.

\bD[bisimulation equivalence]\\
\label{bisimdef}%
Two world-pointed models $\M,w$ and $\M',w'$ are \emph{$n$-bisimilar}, 
$\M,w\,\simn\, \M',w'$, if \PII\ has a winning strategy in the $n$-round 
game starting from $( w,w')$.
$\M,w$ and $\M',w'$ are \emph{bisimilar}, denoted 
$\M,w\sim \M',w'$, if \PII\ has a winning strategy in the unbounded
 game starting from $( w,w')$. 
Two state-pointed models $\M,s$ and $\M',s'$ are ($n$-)bisimilar, denoted $\M,s\!\sim\!\M',s'$
(or $\M,s\,\simn\,\M',s'$), if every world in $s$
is ($n$-)bisimilar to some world in $s'$ and vice versa.  
Two models $\M$ and $\M'$ are \emph{globally bisimilar}, denoted 
$\M\sim \M'$, if every world in $\M$ is bisimilar 
to some world in $\M'$ 
and vice versa.
\eD

\subsection{An inquisitive modal Ehrenfeucht--Fra\"\i ss\'e theorem} 
\label{sec:ef} 

In line with the general situation for back-and-forth games that are
modelled on the expressive power of given logics, the finite levels of
inquisitive bisimulation equivalence capture levels of
indistinguishability in terms of  inquisitive modal logic. This is
brought out in the corresponding Ehrenfeucht--Fra\"\i ss\'e theorem:
for finite sets \P\ of atomic propositions and of agents $\A$, 
{$n$-bisimilarity} $\sim^n$ coincides with $\equiv^n_\inqbm$, 
which is logical indistinguishability in $\inqbm$ at
nesting depth up to~$n$, $\inqbm_n$.
As an immediate consequence, the semantics of all of $\inqbm$ is
preserved under inquisitive bisimulation, which refines all levels
$\sim^n$ (in fact, it even strictly refines their common refinement).

\bT[Ehrenfeucht--Fra\'\i ss\'e theorem for \inqbm\ \cite{CiardelliOttoJSL2021}]\label{EFthm}~\\
\emph{Over any finite sets of agents~$\A$ and atomic propositions \P,
  and for any $n \in \N$ and inquisitive epistemic world- or
  state-pointed models $\M,w$ and $\M',w'$ or $\M,s$ and $\M',s'$:}
 \bre
\item $\M,w\;\simn\; \M',w' \iff
  \M,w \equiv^n_\inqbm\M',w'$
\item $\M,s\;\simn\; \M',s' \iff
\M,s \equiv^n_\inqbm\M',s'$
\ere
\eT

The implications from left to right
go through (by simple syntactic induction) regardless of the finiteness
requirement for $\A$ and $\P$. 
The finiteness requirement is necessary, though,
for the implications from right to left, which can be based on the
provision of characteristic formulae
$\chi^n_{\M,w}$ and $\chi^n_{\M,s} \in \inqbm_n$ that
define the $\sim^n$-equivalence class of $\M,w$, and the
persistent closure of the 
$\sim^n$-equivalence class of $\M,s$, respectively;
\[
 \barr{@{}r@{\;\;\iff\;\;} l@{}}
  \M',w'\models\chi^n_{\M,w}& \M',w'\,\simn\; \M,w
  \\
  \hnt
  \M',s'\models\chi^n_{\M,s} & \M',s'\,\simn\; \M,t
  \mbox { for some } t\subseteq s
\earr
\]
For the details of the proof, the reader is referred to \cite{CiardelliOttoJSL2021}.

Making use of the characteristic formulae, it is easy to obtain the following 
characterisation of the properties definable in \inqbm.

\bC 
\label{EFcorr}
Over a finite signature (i.e.\ with finite sets of agents $\A$, and of propositions $\P$),
a property of world-pointed models
is definable in \inqbm\ if and only if it is closed under $\simn$ for some $n\in\N$. 
Similarly for an inquisitive property of state-pointed models.
\eC

\subsection{Relational inquisitive models}
\label{relinqmodsec}

The main purpose of this paper lies in the characterisation of
inquisitive epistemic logic as a fragment of first-order logic. But
the relevant inquisitive models are not directly conceived as
relational structures (in the manner that e.g.\ ordinary Kripke
structures can directly be thought of as relational structures for
first-order logic). In order to translate the inquisitive setting, 
which involves a universe of information states over and above the
universe of worlds, into a relational structure in which first-order
logic similarly gains access to both kinds of entities, 
we adopt the two-sorted perspective developed in~\cite{CiardelliOttoJSL2021}.
In this view, $W$ and $\wp(W)$ contribute domains for two distinct sorts.
  
\bD[relational models]
\label{def:relational models}
For a signature with sets $\A$ of agents and
$\P= \{ p_i \colon i\in I \}$ of atomic propositions, a 
\emph{relational inquisitive epistemic model} is a two-sorted relational structure
\[
\str{M}=(W,S,(E_a)_{a\in\A},\ee,(P_i)_{i\in I})
\]
with disjoint, non-empty universes
$W,S$ for the two sorts, with binary relations $E_a,\ee\subseteq
W\times S$, and, for $i\in I$, unary $P_i\subseteq W$, satisfying certain conditions
listed below. To state these conditions, we introduce the following notation:
$E_a[w]:= \{ s \in S \colon (w,s) \in E_a\}$ for the set of
$E_a$-successors of $w$; 
$\underline{s} :=\{w\in W \colon (w,s) \in \ee \} \subset W$ for
the set of
$\ee$-predecessors of $s \in
S$; and 
$R_a[w] := \bigcup \{ \underline{s} \colon s \in E_a[w] \}$
for the union of these sets 
across $E_a[w]$ (or, in other words, $R_a:=E_a\circ\ee^{-1}$). 
We can then state the relevant conditions as follows:
\bre
\item extensionality: if $\underline{s}=\underline{s}'$, then $s=s'$.
\item local powerset: if $s\in S$ and $t\subseteq\underline{s}$, there is an $s'\in S$ such that 
  $\underline{s}'=t$.
  \ere
  and for all $a \in \A$ and $w,w' \in W$:
  \bre
  \addtocounter{enumi}{2}
\item non-emptiness: $E_a[w]\neq\emptyset$.
\item downward closure: 
if $s\in E_a[w]$, then $t\in E_a[w]$ for all $t\in S$ with $\underline{t}\subseteq\underline{s}$.
\item factivity:  $w \in  R_a[w]$. 
\item full introspection: if $w' \in R_a[w]$, then $E_a[w'] = E_a[w]$.
\ere
\eD

By extensionality, the second sort $S$ of such a relational model can 
always be identified with a domain of sets over the first sort,
namely, 
$\{\underline{s} \colon s\in S\}
\subset \wp(W)$. 
We always assume this identification and view a 
relational model as a structure $\str{M} =(W,S,(E_a),\in,(P_i))$ where 
$S\subseteq\wp(W)$ and $\in$ is the actual membership relation; this
allows us to specify relational models as just $\str{M}=(W, S, (E_a),(P_i))$.

With a relational model $\str{M}$ we associate the Kripke structure
$\str{K}(\str{M})=(W, (R_a),(P_i)_{i\in I})$,
where $R_a$ is the accessibility relation defined by 
$R_a[w] = \bigcup_{a\in\A} E_a[w]$ (or, equivalently,
by $R_a=E_a\circ\ni$).
Note that by virtue of~(v) and~(vi), each $R_a$ is an equivalence
relation on $W$, and so $\str{K}(\str{M})$ is a (multi-agent) $S5$ model.
The $R_a$-equivalence class of~$w$ is also denoted as $[w]_a  \in W/R_a$. 

\medskip
Further natural constraints on a relational model $\str{M}$ impose
richness conditions on the second sort $S \subset \wp(W)$.
A relational model $\str{M}=(W,S,(E_a),(P_i))$ is said to be
\begin{itemize}
\item[--] \emph{full} if $S=\wp(W)$;
\item[--] \emph{locally full} if for all $a\in\A$ and $w\in W$,
  $\wp([w]_a) \subseteq S$.
\end{itemize}

Any relational model
$\str{M}=(W,S,(E_a),(P_i))$ uniquely determines an inquisitive epistemic
model $\str{M}^\ast=(W,(\Sigma_a),V)$ where $\Sigma_a(w)=E_a[w]$ and
$V(p_i)=P_i$.
This natural passage from relational models to ordinary models induces
semantics for $\inqml$ over relational models and supports related
natural notions like bisimulation equivalence over these.
But as this passage obliterates 
the explicit (extensional) account of
the second sort $S$, different relational models $\str M$
can determine the same inquisitive model $\M$. In other words, 
a given inquisitive epistemic model may come with vastly different
relational counterparts or \emph{relational encodings}. Here a 
\emph{relational encoding} of an inquisitive epistemic model $\M$ is
any relational model $\str M$ such that $\str M^\ast=\M$.

For a given inquisitive epistemic model $\M$, the following three
relational encodings are particularly natural choices, with minimal
or maximal extension for the second sort $S \subset \wp(W)$, and
-- particularly useful for our purposes -- an intermediate version with
minimal extension among locally full encodings.

\bD[relational encodings]
\label{relencdef}
For an inquisitive epistemic model $\M=(W,(\Sigma_a),V)$, we define
the following relational encodings 
of $\M$, each based on~$W$ as the first sort, and always with 
$wE_a s\Leftrightarrow s\in\Sigma_a(w)$, $w\;\ee\; s\Leftrightarrow w\in s$ 
and $P_i\!  =\! V(p_i)$: 
\begin{itemize}
\item[--] $\MM^\rel(\M)$ (minimal) with
$S:= \bigcup_{a\in\A} \mathrm{image}(\Sigma_a)$;
\item[--] $\MM^\lf(\M)$ (minimal locally full) with
  $S:= \bigcup_{a\in\A} \bigcup_{w\in W} \wp([w]_a)$;
\item[--] $\MM^\full(\M)$ (full, maximal) with $S:=\wp(W)$.
\end{itemize}
Corresponding encodings of state-pointed models $(\M,s)$
augment the respective $S$ by $\wp(s)$. 
\eD

To illustrate the difference among the three encodings, imagine that a model $\M$ has universe $W=\{w_1,w_2,w_3\}$ and suppose that for each $a$ and each $i$ we have $\Sigma_a(w_i)=\{\{w_1\},\{w_2\},\emptyset\}$. Then 
\begin{itemize}
\item in the minimal relational encoding $\MM^\rel(\M)$, $S=\{\{w_1\},\{w_2\},\emptyset\}$;
\item in the minimal locally full encoding $\MM^\lf(\M)$, $S=\wp(\{w_1,w_2\})$;
\item in the full encoding $\MM^\full(\M)$, $S=\wp(\{w_1,w_2,w_3\})$.
\end{itemize}

Not surprisingly, the semantics of \inqbm\ over world- or 
state-pointed inquisitive epistemic models translates neatly into 
first-order logic over their two-sorted relational encodings, 
uniformly and regardless of the specific choices for the second sort. 
A systematic standard translation, as specified in the
general inquisitive modal setting in~\cite{CiardelliOttoJSL2021},
follows the defining clauses in Definition~\ref{supportsemdefn},
reformulated in terms of first-order quantification over the
first or second sort of a relational encoding, respectively.%
\footnote{There is a small subtlety here. In \cite{CiardelliOttoJSL2021},
two translations are introduced. The first directly mirrors the clauses of 
the semantics, but in general it only works on encodings which are locally fully. 
The second, $\phi\mapsto\phi^\star$,
works on arbitrary models, but it requires 
a pre-processing of the formula $\phi$ into an equivalent $\Box$-free formula. 
We refer to ~\cite{CiardelliOttoJSL2021} for the details.} 
Denoting
this \emph{standard translation} by $\phi \mapsto \phi_{\texttt{w}}^\star$
(with first-order variable symbol $\texttt w$ for worlds, over the first sort) or 
$\phi \mapsto \phi_{\texttt{s}}^\star$
(with first-order variable symbol $\texttt s$ for information states, over the second sort),
we find the following.

\bO\label{obs:correspondence full}
For any inquisitive epistemic model $\M$, any relational 
encoding $\str{M}$ of $\M$, the first-order standard translations
$\phi_{\emph{\texttt{w}}}^\star$ 
and
$\phi_{\emph{\texttt{s}}}^\star$ 
of 
$\phi\in\inqml$ are such that for all worlds $w\in W$ and all states $s\in S$:
\bre
\item $\M,w\models\phi\iff\MM\models\phi_{\emph{\texttt{w}}}^\star[w]$
\item $\M,s\models\phi \iff \MM\models\phi_{\emph{\texttt{s}}}^\star[s]$
\ere
\eO

Modulo this natural standard translation, $\inqbm$ can be viewed as a
syntactic fragment of first-order logic,
\[
  \inqbm \subset \FO,
\]
over the class of all relational inquisitive models. This situation
is similar to what the classical standard translation achieves
for modal logic $\ML$ over Kripke structures, which are naturally
conceived as relational structures. The crucial difference for the
inquisitive setting is that the target class of structures,
unlike the class of relational Kripke structures, forms a
non-elementary class. The non-elementary nature of any natural
class of relational inquisitive encodings stems form the
downward closure condition for inquisitive assignments. This requires
representation of the full power set of 
each information state in the image of $\Sigma$ in the second sort.
That this constraint cannot be captured by $\FO$ conditions, can technically
be demonstrated by means of compactness arguments for $\FO$.
We refer to~\cite{CiardelliOttoJSL2021} for related examples,
and to~\cite{meissnerOtto22} for a model-theoretic account that
relaxes the notion of (relational) models via passage to a wider class
of (relational) pseudo-models, which is elementary and, in a precise sense,
conservative w.r.t.\ \inqbm\ semantics.

As a fragment of first-order logic over relational models, i.e.\ over
relational encodings of the intended models, 
$\inqbm \subset \FO$ displays crucial model-theoretic preservation
properties, viz.\ 
$\sim$-invariance and, in the state-pointed case, 
persistency and the empty state property. 
In showing that $\inqbm \subset \FO$  is
\emph{expressively complete} for all those $\FO$-properties of worlds
or states (in relational inquisitive epistemic models from suitable background
classes $\CC$) that satisfy these semantic
constraints, we obtain the desired \emph{characterisation} result.
In other words, we want to show that if a $\FO$-formula
$\phi(\texttt w)$ in a single free variable of the first sort
 defines a property of worlds which is bisimulation invariant over \CC,
 then it is logically equivalent over $\CC$ to some formula of $\inqbm$
(in defining a world property);  and similarly,
that if a $\FO$-formula $\phi(\texttt s)$ in a single free variable of the second 
sort defines a property of information states which, over~$\CC$,
is bisimulation invariant and 
inquisitive (in the sense of Definition \ref{inqclass}), 
then it is logically equivalent over $\CC$ to some formula of $\inqbm$
(in defining a state property). In a nutshell, we seek to establish the semantic
equality
\[
\inqbm \equiv \FO/{\sim} \quad \mbox{ $(\dagger)$}
\]
over suitable classes $\CC$ of relational epistemic models. Each of these
equivalences in expressiveness comes in two readings, one for world
properties and one for state properties. We shall mostly in the first
instance emphasise the simpler world-pointed reading and discuss the
state-pointed versions second. 

Such expressive completeness
and characterisation results, for world properties in particular,
match van~Benthem's classical characterisation of modal logic as the
bisimulation-invariant fragment of first-order logic.
Versions of this characterisation for basic inquisitive modal logic
have been expounded in~\cite{CiardelliOttoJSL2021}; those results
establish the desired expressive completeness results over the classes
of (i) all relational inquisitive models, (ii) all finite relational inquisitive models
(iii) all locally full models, and (iv) all finite locally full models.%
\footnote{Relational inquisitive models are defined as in Definition \ref{relencdef},
but without the requirements of factivity and introspection.}
A similar expressive completeness result over the class of all full
relational inquisitive models can be ruled out by a compactness
argument~\cite{CiardelliOttoJSL2021}. 
Due to the non-elementary nature of the relevant classes of
two-sorted relational models, the positive results of~\cite{CiardelliOttoJSL2021}
already need to employ non-classical model-theoretic arguments.
Our present analysis for the epistemic case, however, requires a more
complex route to expressive completeness, due to the additional
constraints imposed on inquisitive epistemic frames, viz. factivity and full introspection. 

In preparation for those arguments, we review some salient technical
background on the model-theoretic analysis of first-order
expressiveness over the kind of two-sorted structures involved here. 
Important tools in this are Ehrenfeucht-Fra\"\i ss\'e techniques
for first-order logic $\FO$ as well as for monadic second-order logic $\MSO$. 

\subsection{Classical Ehrenfeucht--Fra\"\i ss\'e techniques}
\label{EFtechniques}

We review the classical Ehrenfeucht--Fra\"\i ss\'e analysis of the expressive
power of $\FO$ over relational models. As our relational models are two-sorted, with a
second sort $S \subset \wp(W)$ that represents subsets of the first
sort $W$, we shall have occasion to look beyond plain first-order expressiveness
and to consider also limited forms of monadic second-order
expressiveness. The limitations that account for the tame nature of
$\MSO$ applications in our analysis stem from their effective restriction to
the local subuniverses $[w]_a$. 
We note that, for locally full relational models, first-order
quantification over the second sort can indeed emulate full $\MSO$
in restriction to $[w]_a \subset W$, since for these models 
$\wp([w]_a) \subset S$.

\paragraph*{Ehrenfeucht--Fra\"\i ss\'e for FO.}
We review the classical Ehrenfeucht--Fra\"\i ss\'e analysis of the expressive
power of $\FO$ over relational structures. The underlying game, like
the bisimulation game, is a two-player game, that serves to probe, in this
case, levels of equivalence w.r.t.\ first-order logic for matched tuples
of elements over these two structures. Positions in the game over
relational structures $\str{A},\str{A}'$ are pairs $(\abar;\abar')$ with
marked tuples $\abar$ of elements in $\str{A}$ and $\abar'$ of elements
in $\str{A}'$ of the same finite length; 
we may think of $\abar$ and $\abar'$
as tuples of elements that are pebbled with pairs of matching pebbles so
as to match the components of $\abar$ to the components of $\abar'$.
In a single round, player~\PI\ positions one pebble of the next pair of
pebbles in either $\str{A}$ or $\str{A}'$, and player~\PII\ must position
its match in the opposite structure. This takes the game from position
$(\abar;\abar')$ to some successor position $(\abar,a;\abar',a')$.
Player~\PII\ loses in position  $(\abar;\abar')$ if
$\str{A}\restr \abar, \abar \not\simeq \str{A}'\restr \abar', \abar'$,
i.e.\ if the pebbled match does not define an isomorphism between 
the induced substructures over these elements. 
If no such position is reached by the end of the game, \PII\ wins.
A winning  strategy
for \PII\ in the $n$-round game starting from position $(\abar;\abar')$ over
$\str{A}$ and $\str{A}'$  guarantees that
$\str{A},\abar \equiv_n \str{A}',\abar'$ ($n$-elementary equivalence),
which means that $\abar$ in $\str{A}$ and $\abar'$ in $\str{A}'$ are
indistinguishable by $\FO$-formulae of quantifier rank up to $n$.
Denoting the equivalence relation defined in terms of a
winning  strategy for \PII\ in the $n$-round game by $\simeq_n$:
\begin{equation}
\str{A},\abar \simeq_n \str{A}',\abar'
\quad \Rightarrow \quad \str{A},\abar \equiv_n \str{A}',\abar'
\label{EFhalfFOeqn}
\end{equation}

For structures $\str{A}$ and $\str{A}'$ in a finite relational signature,
the classical Ehren\-feucht--Fra\"\i\-ss\'e theorem states that the
converse implication also holds, so that~(\ref{EFhalfFOeqn}) becomes an
equivalence, which characterises levels of elementary equivalence in
terms of the associated back-and-forth game.

\medskip
The precise correspondence in the Ehrenfeucht--Fra\"\i ss\'e theorem witnesses how
the positioning of a next pebble pair in a round of the
game for $\FO$ probes the behaviour of the currently matched tuples
w.r.t.\ a single further $\forall/\exists$-quantification over
elements --
corresponding to an increase of~$1$ in terms of quantifier rank. The
corresponding game for $\MSO$ will additionally allow \PI\ to probe
equivalence w.r.t.\
$\forall/\exists$-quantification over subsets, and game positions
correspondingly record a match between elements \emph{and} subsets.
In our specific context we shall exploit parallels beween the
$\FO$-game involving both sorts of our relational
models and the $\MSO$-game over the first sort.

\paragraph*{Ehrenfeucht--Fra\"\i ss\'e for MSO.}
We consider two ordinary structures
$\str{A},\str{A}'$ in the same relational signature, with
(single-sorted) universes $A,A'$.
To capture the expressiveness of $\MSO$ over these, game
positions now consist of pairings $(\Pbar,\abar;\Pbar',\abar')$ where
$\abar$ and $\abar'$ are tuples of elements of the same finite length
from $A$ and $A'$, and $\Pbar$ and $\Pbar'$ are tuples of subsets 
of $A$ and $A'$ of the same finite length. The condition that \PII\
must maintain in order not to lose in position 
$(\Pbar,\abar;\Pbar',\abar')$ over $\str{A},\str{A}'$
now is 
\[
  (\str{A},\Pbar)\restr\abar,\abar \simeq (\str{A}',\Pbar')\restr\abar',\abar'.
\]
This means that the elements matched by the pebble pairs
mark out isomorphic induced substructures in the expansions of
$\str{A}$ by the unary predicates $\Pbar$  and of $\str{A}'$ by $\Pbar'$.
These expansions may be seen as colourings produced in second-order
moves in which $\Pbar$ and $\Pbar'$ were matched, which are probed
(just as in the first-order $\FO$ game) by the match in pebbled
elements $\abar$ and $\abar'$ produced in first-order moves. In each
round, \PI\ determines whether to play it in first- or second-order
mode, and correspondingly makes a choice of either an element or a
subset on the $\str{A}$-side or  on the $\str{A}'$-side; \PII\ must
respond likewise on the opposite side. Depending on the mode of this
round, it may overall lead from $(\Pbar,\abar;\Pbar',\abar')$ to
$(\Pbar,\abar,a;\Pbar',\abar',a')$ (with a new pair of pebbled
elements $(a,a')$ for first-order round)
or to $(\Pbar,P,\abar;\Pbar',P',\abar')$ (with a new pair of subsets or
a new colour for elements, $(P,P')$ for a second-order round).

A winning  strategy for \PII\ for $n$ rounds of this game
starting from position $(\Pbar,\abar;\Pbar',\abar')$ over
$\str{A}$ and $\str{A}'$ now guarantees that
$\str{A},\Pbar,\abar \equiv_n^\MSO \str{A}',\Pbar',\abar'$, i.e.\
indistinguishability by $\MSO$-formulae $\phi(\Xbar,\xbar)$ of 
quantifier rank up to $n$.%
\footnote{We here appeal to a `mixed' quantifier rank for $\MSO$ that does not
  discriminate between first- and second-order quantification.}
Denoting the equivalence relation defined in terms of a
winning  strategy for \PII\ in this $n$-round game for $\MSO$ by
$\simeq_n^\MSO$, the general result analogous to~(\ref{EFhalfFOeqn}) above is
\begin{equation}
\str{M},\mbar \simeq_n^\MSO \str{M}',\mbar'
\quad \Rightarrow \quad \str{M},\mbar \equiv_n^\MSO \str{M}',\mbar' 
\label{EFhalfMSOeqn}
\end{equation}
and, again, the converse implication is
also true over finite relational signatures.

\section{\texorpdfstring{\boldmath
  Expressive completeness for epistemic $\inqbm$}{Expressive completeness for epistemic InqML}}
\label{epinqsec}

\subsection{Bisimulation invariance as a semantic constraint} 
\label{bisimconstraintsec}

After the introduction of all the relevant notions, our main
theorem announced in Section~\ref{introsec} can be stated more
precisely as follows. 

\medskip\noindent
\textsc{Theorem~\ref{main2}}. \label{th:restated} 
\emph{Let $\CC$ be any one of the following classes of relational models: 
the class of all inquisitive epistemic models; of all
finite inquisitive epistemic models; of all locally full,
or of all finite locally full inquisitive epistemic models.
Over each of these classes,
\[
\inqbm \equiv \FO/{\sim} \quad \mbox{ $(\dagger)$}
\]
i.e., a property of world-pointed models is definable
in \inqbm\ over $\CC$ if and only if it is both $\FO$-definable over $\CC$ 
and $\sim$-invariant over  $\CC$.%
\footnote{As usual, $\FO$-definability of a property $\P$ of world-pointed
  models over a class $\CC$ of models means that there is a
  $\FO$-formula $\phi$ with a single free variable of the first sort
such that for all $\MM \in \CC$:  $(\MM,w) \in \P$ iff 
$\MM\models\phi(w)$; 
similarly for $\inqbm$-definability over $\CC$ and for properties of
state-pointed models and formulae with a free variable of the second sort.} 
Similarly, an inquisitive property of state-pointed models is definable in \inqbm\ over $\CC$ if and only if it is $\FO$-definable over $\CC$ and $\sim$-invariant over $\CC$.}

\medskip
The inquisitive Ehrenfeucht--Fra\"\i ss\'e theorem, Theorem~\ref{EFthm}, 
implies that $\inqbm$ and its standard translation into $\FO$ 
are invariant under bisimulation: $\inqbm \subset \FO/{\sim}$.
By Corollary~\ref{EFcorr} it further implies \emph{expressive completeness}
of $\inqbm_n$ for any $\simn$-invariant property of world-pointed
models and any inquisitive $\simn$-invariant property of
state-pointed models: $\FO/{\simn} \subset \inqbm$ for every $n \in \N$. 
In order to prove $(\dagger)$ in restriction to some particular class $\CC$ of
relational inquisitive models in the world-pointed case, 
it therefore remains to show that, for any formula 
$\phi(\w) \in \FO$ in a single free variable of the first sort, there
is some $n \in \N$ such that $\sim$-invariance of $\phi(\w)$ over $\CC$ implies 
 $\simn$-invariance over $\CC$.
 \[
   \phi(\w) \in \FO/{\sim} \; \Longrightarrow \; \phi(\w) \in \FO/{\simn} \mbox{
     for some $n = n(\phi)$ } \quad 
   \mbox{ $(\ddagger)$}
\]

This may be viewed as a non-trivial \emph{compactness principle} for 
$\sim$-invariance of first-order properties, over the non-elementary
classes $\CC$ of interest.

\bO
\label{compactobs}
\bae
\item
For any class $\CC$ of 
relational inquisitive models, the following are equivalent:
\bre
\item
$\inqbm \equiv \FO/{\sim}$ for world properties over $\CC$;
\item
for $\FO$-properties of world-pointed models, 
$\sim$-invariance over $\CC$ implies $\simn$-invariance over \CC\ 
for some~$n$.
\ere
\item
Similarly, the following are equivalent:
\bre
\item
$\inqbm \equiv \FO/{\sim}$ for 
inquisitive state properties over~$\CC$;
\item
for inquisitive
$\FO$-properties of state-pointed models, 
$\sim$-in\-vari\-ance over $\CC$ implies $\simn$-invariance over \CC\ 
for some~$n$.
\ere
\eae
\eO

As noted above, examples
in~\cite{CiardelliOttoJSL2021} illustrate that $\FO$ does not satisfy 
a compactness theorem over the non-elementary class of relational
inquisitive models. Moreover, over the class of 
\emph{full} inquisitive relational models, failures of compactness can be
exhibited even for the bisimulation-invariant fragment of $\FO$,
which must therefore be a strict extension of $\inqbm$ over this class.
In other words, an analogue of Theorem~\ref{main2} fails for the class of full models.

\begin{figure}
\[
\xymatrix{
\str{M},w \ar@{-}[d]|{\rule{0pt}{1ex}\sim} \ar @{-}[rr]|{\;\simn\,}
&& \str{M}'\!,w' \ar @{-}[d]|{\rule{0pt}{1ex}\sim}
\\
\hat{\str{M}},\hat{w} \ar@{-}[rr]|{\rule{0pt}{1.5ex}\;\equiv_q}
&& \hat{\str{M}}',\hat{w}'
}
\qquad\qquad
\xymatrix{
\str{M},s \ar@{-}[d]|{\rule{0pt}{1ex}\sim} \ar @{-}[rr]|{\;\simn\,}
&& \str{M}'\!,s' \ar @{-}[d]|{\rule{0pt}{1ex}\sim}
\\
\hat{\str{M}},\hat{s} \ar@{-}[rr]|{\rule{0pt}{1.5ex}\;\equiv_q}
&& \hat{\str{M}}'\!,\hat{s}'
}
\]
\caption{Generic upgrading patterns.}
\label{genupgradefigure}
\end{figure}
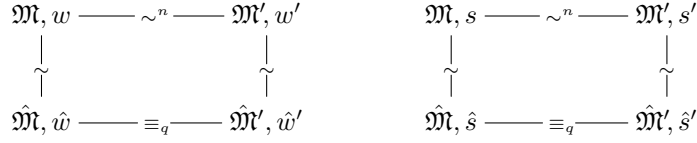

To establish~$(\dagger)$ 
over some class $\CC$ of relational models, we aim to show
condition~$(\ddagger)$ by means of an upgrading of $\simn$-equivalence
to $\equiv_q$, where $q$ is the quantifier rank of the given formula
$\phi \in \FO$ and $n=n(q)$ is a suitable natural number depending on $q$. 
If suitable bisimilar partner structures like
$\hat{\str{M}}$ and $\hat{\str{M}}'$ in Figure~\ref{genupgradefigure}
are available within $\CC$ for any pair of $\simn$-equivalent pointed
structures $\str{M}$ and $\str{M}'$ from $\CC$, then the detour
through the lower rung in these diagrams shows the semantics of
$\phi$ to be $\simn$-invariant over $\CC$, hence equivalently
expressible in $\inqbm$ over $\CC$ due to Observation~\ref{compactobs}. 

This general approach to expressive completeness has been tested
in several cognate scenarios that range from a simple and constructive
proof in~\cite{OttoNote,GorankoOtto}
of van~Benthem's characterisation of $\ML$~\cite{Benthem83} and
Rosen's finite model theory version~\cite{Rosen}
to the characterisation of $\inqbm$ as the
$\sim$-invariant fragment of $\FO$ over (finite or general, possibly
just locally full) relational inquisitive models in the basic rather than the epistemic
setting in~\cite{CiardelliOttoJSL2021}.

The upgrading argument needs to apply within $\CC$, where
$\sim$-invariance can be assumed for the vertical
transfers. At the same time it must overcome,
with the choice of $\hat{\str{M}}$ and $\hat{\str{M}}'$,
obstacles to the desired $\FO$-indistinguishability that are not
governed by $n$-bisimilarity. Among those features are issues related
to distinctions in low-level multiplicities or local cyclic patterns
that can be expressed in $\FO$ at quantifier rank~$q$. 
So within the classes $\CC$ of relational epistemic models that are of
interest here, we need to come up with bisimilar partner structures $\hat{\str{M}}$ and
$\hat{\str{M}}'$, for which Ehrenfeucht--Fra\"\i ss\'e arguments 
can establish $\equiv_q$.
That task can be further split into two main stages:
\bne
\item
  the construction of suitable bisimilar inquisitive epistemic
  companion models $\hat{\str{M}}\sim \str{M}$ and $\hat{\str{M}}' \sim \str{M}'$,
and   
\item
  a model-theoretic interpretation of the relevant relational
  encodings
$\hat{\str{M}}$ and $\hat{\str{M}}'$
of $\hat{\M}$ and $\hat{\M}'$ in related two-sorted first-order structures
$\expl{\hat{\M}} = \expl{\hat{\str{M}}}$ and $\expl{\hat{\M}'} = \expl{\hat{\str{M}}'}$;
  this yields a neater argument for $\equiv_q$ involving, among other
  techniques, locality arguments based on Gaifman's theorem for~$\FO$.
\ene

\begin{figure}
\[
\xymatrix{
\str{M},w \ar@{-}[d]|{\rule{0pt}{1ex}\sim}\ar @{-}[rr]|{\;\simn\,}
&& \str{M}'\!,w' \ar @{-}[d]|{\rule{0pt}{1ex}\sim}
\\
\hat{\str{M}},\hat{w} \ar@{-}[rr]|{\rule{0pt}{1.5ex}\;\equiv_q}
\ar@{<->}[d]_{\expl{\cdot}}
&& \hat{\str{M}}',\hat{w}'
\ar@{<->}[d]^{\expl{\cdot}}
\\
\expl{\hat{\str{M}},\hat{w}} \ar@{-}[rr]|{\rule{0pt}{1.5ex}\;\equiv_{q+d}}
&& \expl{\hat{\str{M}}',\hat{w}'}
}
\qquad\qquad
\xymatrix{
\str{M},s \ar@{-}[d]|{\rule{0pt}{1ex}\sim} \ar @{-}[rr]|{\;\simn\,}
&& \str{M}'\!,s' \ar @{-}[d]|{\rule{0pt}{1ex}\sim}
\\
\hat{\str{M}},\hat{s} \ar@{-}[rr]|{\rule{0pt}{1.5ex}\;\equiv_q}
\ar@{<->}[d]_{\expl{\cdot}}
&& \hat{\str{M}}',\hat{s}'
\ar@{<->}[d]^{\expl{\cdot}}
\\
\expl{\hat{\str{M}},\hat{s}} \ar@{-}[rr]|{\rule{0pt}{1.5ex}\;\equiv_{q+d}}
&& \expl{\hat{\str{M}}',\hat{s}'}
}
\]
\caption{Upgrading patterns for relational epistemic models.}
\label{S5upgradefigure}
\end{figure}
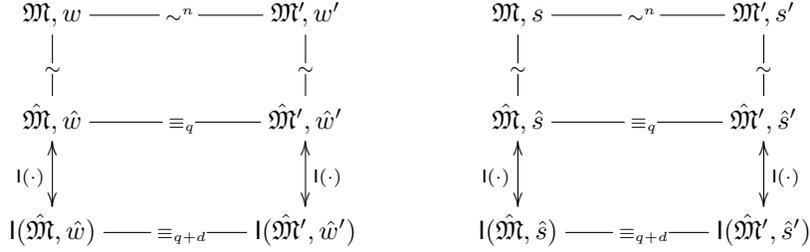

We first review and prepare the ground for~(2) in Section~\ref{FOEFGaifsec}
below; then address~(1) in Section~\ref{FOscattersec}; and finally combine
these preparatory steps to prove our main characterisation  theorem,
Theorem~\ref{main2}, in Section~\ref {epicharsec}.
The picture for the intended upgrading arguments can be 
extended as shown in Figure~\ref{S5upgradefigure}.

\subsection{First-order expressiveness over relational models} 
\label{FOEFGaifsec}

As pointed out above, the assessment of the expressive power of 
$\FO$ over relevant non-elementary classes of relational structures
cannot use classical compactness arguments. Instead it draws on
Ehrenfeucht--Fra\"\i ss\'e arguments. These rely on 
\emph{locality criteria} based on \emph{Gaifman distance}
(cf.~Definition~\ref{Gaifmandef}) and suitable local tree-like
unfoldings in appropriate image structures $\expl{\str{M}}$ --
structural correlates that lend themselves 
to the localisation of Ehrenfeucht--Fra\"\i ss\'e arguments, as it
were, and related to $\str{M}$ by mutual $\FO$-interpretability.

\paragraph*{First-order interpretations in an exploded view.}

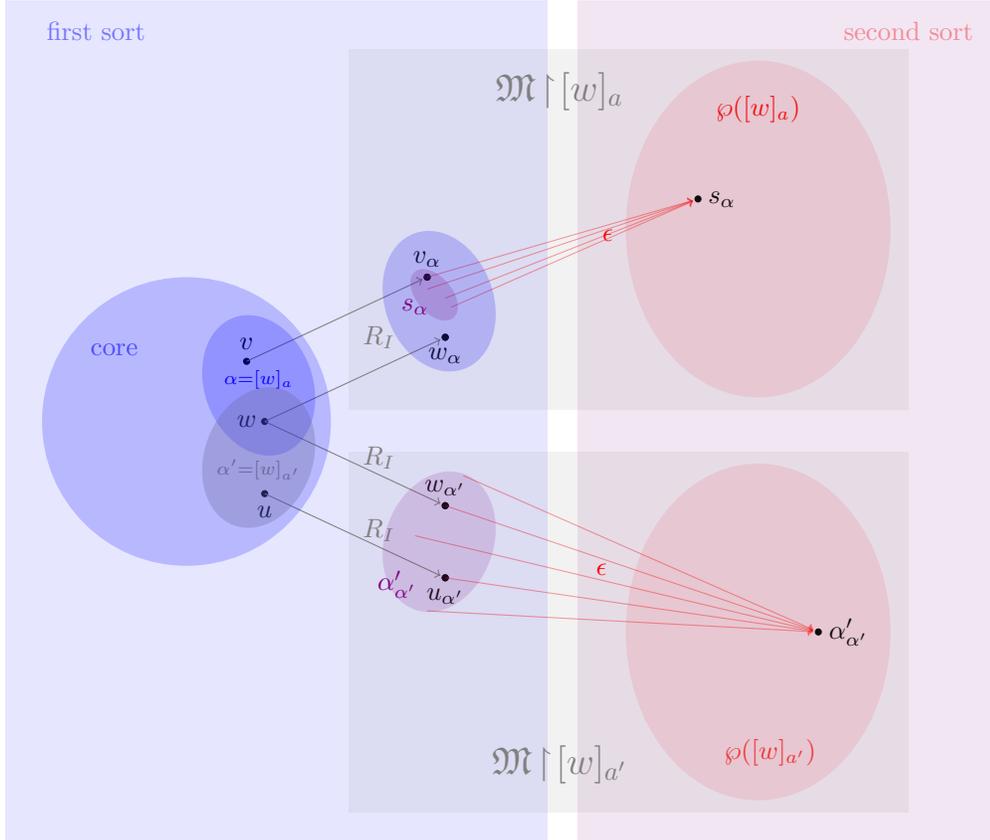
\begin{figure}
\[
{\begin{tikzpicture}[scale=.8]
\node at (-7.5,11.5) [color=blue, opacity=.5]  {first sort};
\node at (6,11.5) [color=red, opacity=.4]  {second sort};
\path[fill=blue, fill opacity=.1] (-9,-2) rectangle ++(9cm,14cm);
\path[fill=violet, fill opacity=.1] (.5,-2) rectangle ++(7cm,14cm);
\coordinate (core) at (-6,5);
\node at ($(core)+(-1.2,1.2)$) [color=blue, opacity=.6]  {core};
\coordinate (a1) at ($(core)+(1.2,.6)$);
\coordinate (a2) at ($(core)+(1.2,-.6)$);
\coordinate (c1) at ($(a1)+ (3,1.4)$);
\coordinate (c2) at ($(a2)+(3,-1.4)$);
\coordinate (p1) at ($(core)+(9.5,3.2)$);
\coordinate (sa1) at ($(p1)+(-1,.5)$);
\coordinate (p2) at ($(core)+(9.5,-3.5)$);
\coordinate (sa2) at ($(p2)+(1,0)$);
\coordinate (w) at ($(core)+(1.3,0)$);
\coordinate (v) at ($(core)+(1,1)$);
\coordinate (u) at ($(core)+(1.3,-1.2)$);
\coordinate (wc1) at ($(w)+(3,1.4)$);
\coordinate (wc2) at ($(w)+(3,-1.4)$);
\coordinate (vc1) at ($(v)+(3,1.4)$);
\coordinate (uc2) at ($(u)+(3,-1.4)$);
\coordinate (Malpha1) at ($(c1)+ (2.5,.5)$);
\coordinate (Malpha2) at ($(c2)+(2.5,-.7)$);
\node at ($(Malpha1)+(-.5,3)$)  [gray] {\Large $\str{M}\restr [w]_a$};
\node at ($(Malpha2)+(-.5,-3)$) [gray] {\Large $\str{M}\restr [w]_{a'}$};
\node at ($(w)+(-.3,0)$) {$w$};
\node at ($(v)+(0,.3)$) {$v$};
\node at ($(u)+(0,-.3)$) {$u$};
\node at ($(wc1)+(0,-.3)$) {$w_\alpha$};
\node at ($(wc2)+(0,.3)$) {$w_{\alpha'}$};
\node at ($(vc1)+(0,.3)$) {$v_\alpha$};
\node at ($(uc2)+(0,-.3)$) {$u_{\alpha'}$};
\draw [fill=black] (u) circle (.05cm);
\draw [fill=black] (v) circle (.05cm);
\draw [fill=black] (w) circle (.05cm);
\draw [fill=black] (uc2) circle (.05cm);
\draw [fill=black] (vc1) circle (.05cm);
\draw [fill=black] (wc1) circle (.05cm);
\draw [fill=black] (wc2) circle (.05cm);
\draw [fill=black] (uc2) circle (.05cm);
\draw [fill=black] (vc1) circle (.05cm);
\draw [fill=black] (wc1) circle (.05cm);
\draw [fill=black] (wc2) circle (.05cm);
\draw [gray][-{>[sep=2pt]}] (v) -- (vc1);
\draw [gray] [-{>[sep=2pt]}] (w) -- (wc1);
\draw [gray] [-{>[sep=2pt]}] (w) -- (wc2);
\draw [gray] [-{>[sep=2pt]}] (u) -- (uc2);
\draw [red,opacity=.4][-{>[sep=2pt]}] (vc1) -- (sa1);
\draw [red,opacity=.4][-{>[sep=2pt]}] ($(vc1)+(0,-.2)$) -- (sa1);
\draw [red,opacity=.4][-{>[sep=2pt]}] ($(vc1)+(.3,-.35)$) -- (sa1);
\draw [red,opacity=.4][-{>[sep=2pt]}] ($(vc1)+(.4,-.5)$) -- (sa1);
\draw [red,opacity=.4][-{>[sep=2pt]}] (uc2) -- (sa2);
\draw [red,opacity=.4][-{>[sep=2pt]}] (wc2) -- (sa2);
\draw [red,opacity=.4][-{>[sep=2pt]}] ($(wc2)+(-.5,-.5)$) -- (sa2);
\draw [red,opacity=.4][-{>[sep=2pt]}] ($(wc2)+(.3,.5)$) -- (sa2);
\draw [red,opacity=.4][-{>[sep=2pt]}] ($(c2)+(-.2,-1.15)$) -- (sa2);
\node at ($(a1)+(2,.8)$) [gray] {$R_I$};
\node at ($(a2)+(2,-1.2)$) [gray] {$R_I$};
\node at ($(a2)+(2,0)$) [gray] {$R_I$};
\def\elength{1.2cm}
\def\ewidth{.9cm}
\path[fill=blue, fill opacity=0.2, rotate=290] (a1) ellipse
({\elength} and {\ewidth});
\node at ($(a1)+(0,.1)$) [blue] {$\scc \alpha= [w]_a$};
\path[fill=gray, fill opacity=0.3, rotate=70] (a2) ellipse ({\elength} and {\ewidth});
\node at ($(a2)+(0,-.2)$) [gray] {$\scc \alpha' = [w]_{a'}$};
\def\radius{2.4cm}
\path[fill=blue, fill opacity=0.2] (core) ellipse ({\radius} and
{\radius});
\path[fill=blue, fill opacity=0.2, rotate=290] (c1) ellipse ({\elength} and {\ewidth});
\path[fill=violet, fill opacity=0.15, rotate=70] (c2) ellipse ({\elength} and {\ewidth});
\def\elength{2.8cm}
\def\ewidth{2.2cm}
\path[fill=red, fill opacity=0.1, rotate=90] (p1) ellipse ({\elength}
and {\ewidth});
\path[fill=red, fill opacity=0.1, rotate=90] (p2) ellipse ({\elength}
and {\ewidth});
\node at ($(p1)+(0,2)$) [red] {$\wp([w]_a)$};
\draw [fill=black] (sa1) circle (.05cm);
\node at ($(sa1)+(.4,0)$) {$s_\alpha$};
\draw [fill=black] (sa2) circle (.05cm);
\node at ($(p2)+(.2,-2)$)[red, opacity=.8] {$\wp([w]_{a'})$};
\node at ($(sa2)+(.5,0)$) {$\alpha'_{\alpha'}$};
\def\elength{.5cm}
\def\ewidth{.3cm}
\path[fill=violet, fill opacity=0.2, rotate=310] ($(vc1)+(.3,-.1)$) ellipse
({\elength} and {\ewidth});
\path[fill=gray, fill opacity=0.1 ] ($(Malpha1)+(-4,-2.3)$) rectangle ++(9.3cm,6cm);
\path[fill=gray, fill opacity=0.1 ] ($(Malpha2)+(-4,-3.8)$) rectangle ++(9.3cm,6cm);
\node at ($(vc1)+(-.2,-.5)$) [violet] {$s_\alpha$};
\node at ($(vc1)+(3,.7)$) [red] {$\epsilon$};
\node at ($(c2)+(-.7,-.7)$) [violet] {$\alpha'_{\alpha'}$};
\node at ($(uc2)+(2.6,0.15)$)[red]{$\epsilon$};
\end{tikzpicture}}
\]
\caption{Structural layout of the exploded view $\expl{\str{M}}$ of
  $\str{M}$; selectively displayed are 
  the representations of intersecting
  local structures $\str{M}\restr [w]_a$, $\str{M}\restr [w]_{a'}$ w.r.t.\ agents $a$
  and $a'$,  with one element each from their (full powerset) second
  sorts.
  Here $\alpha = [w]_a$ and $\alpha'= [w]_{a'}$ contribute two
  overlapping local structures, whose disjoint (tagged)
  representations each contribute the full power set of their first sort
  to disjoint (tagged) representations in the second sort; displayed
  as particular elements of the second sort are $s_\alpha$ for some
  $s \subset  \alpha$ (as part of the representation of the $a$-local
  structure $\str{M}\restr[w]_a$)
  and the full set $\alpha' \subset \alpha'$ (as part of the
  representation of the $a'$-local structure $\str{M}\restr[w]_{a'}$).}
\label{explfigure}
\end{figure}

Consider an epistemic model  
$\M = (W,(\Sigma_a)_{a\in\A},V)$ with locally full two-sorted relational encoding
 $\str{M} := \str{M}^{\lf}(\M) = (W,S, (E_a)_{a \in \A},(P_i)_{i\in I})$, where 
 \[
   S = \bigcup \bigl\{ \wp([w]_a) \colon w\in W, a\in\A \bigr\},
 \]
where $[w]_a = \sigma_a(w) = R_a[w]$.  
 We also consider the natural restrictions-cum-reducts $\M\restr [w]_a$ and 
 $\str{M}\restr[w]_a$ for any combination of $w \in W$ and $a \in \A$:
 \[
     \M\restr[w]_a := \bigl( [w]_a, \Sigma_a \restr[w]_a, V \restr[w]_a \bigr)
 \]
will be referred to as an \emph{$a$-local structure} or just
\emph{local structure} and is treated as inquisitive model for the
single agent~$a$ with constant $\Sigma_a$.
The relational encoding
\[
  \str{M}\restr[w]_a := \bigl( [w]_a, \wp([w]_a), E_a \restr[w]_a,
  (P_i \restr[w]_a) \bigr)
\]
of this local structure
arises as the induced substructure of the reduct of $\str{M}$
that eliminates all other agents; this restriction is two-sorted and full, as the
domain for the second sort, the restriction of $S$ to 
$\wp([w]_a)$ comprises all of $\wp([w]_a)$ since $\str{M} = \str{M}^\lf(\M)$
is locally full. So each $\str{M}\restr[w]_a$ is a full relational encoding of the 
single-agent, single-equivalence-class inquisitive epistemic model $\M\restr [w]_a$.
In order to analyse the expressive power of first-order logic over 
$\str{M}$, we look at an ``exploded view'' $\expl{\str{M}}$ of $\str{M}$,
which represents $\str{M}$ as the conglomerate of these restrictions
$\str{M}\restr [w]_a$, with their interconnections charted in a core
structure representing the underlying Kripke frame. 
The concrete format of such an exploded view
would admit quite some variations. What matters for us is
first-order interpretability of $\str{M}$ itself within
$\expl{\str{M}}$: as indicated in Figure~\ref{S5upgradefigure},
we want to apply these representations to suitably prepared
models $\hat{\str{M}}$ so as to argue for $\equiv_q$ at the level
of $\hat{\str{M}}/ \hat{\str{M}}'$ on the basis of $\equiv_{q+d}$
at the level of the $\expl{\hat{\str{M}}}/\expl{\hat{\str{M}}'}$ for some
fixed constant $d$. We first discuss suitable representations
$\expl{\str{M}}$ for arbitrary locally full relational epistemic
models of the form $\str{M} = \str{M}^\lf(\M)$ but note that they will
only be used for specifically prepared bisimilar partner structures as
in Figure~\ref{S5upgradefigure}; it is also important to note that the
$\expl{\cdot}$-images are not themselves relational encodings of 
inquisitive models.
Figure~\ref{explfigure} illustrates the structural ingredients in the format that
we choose for this exploded view $\expl{\str{M}}$ of $\str{M}$.
The two-sorted structure $\expl{\str{M}}$, from which $\str{M}$ itself
can be retrieved by first-order means,
involves (see below for details):
\bre
\item[--]
a \emph{core part} within the first sort, which represents the underlying
multi-agent $S5$ Kripke frame on the set of worlds $W$, with relation
symbols $R_a$ for accessibility w.r.t.\ agent~$a \in \A$ as in the
induced Kripke frame;
\item[--]
a halo of attached disjoint external copies of
each one of the two-sorted \emph{local structures} $\str{M}\restr [w]_a$, one for each
equivalence class of worlds w.r.t.\ to agents $a \in \A$; $\str{M}\restr [w]_a$
carries the natural restrictions of the relations in $\str{M}$ but with inquisitive
assignment just for agent~$a$ (see above); 
\item[--]
with additional unary marker predicates~$P_a$ for the $a$-local
components and a binary  relation~$R_I$ for the identification
between worlds in the core part and their instances in local
component structures in the halo.%
\ere

The first part of the disjoint union underlying $\expl{\str{M}}$
forms the core (fully in the first sort), the remaining parts (all truly two-sorted)
form the halo, with the interpretation of $R_I$ as the only link between
the otherwise disjoint pieces:
\[
    \expl{\str{M}} := 
    \displaystyle
    \bigl(
    (W, (R_a)_{a \in \A}) 
\;\oplus \!\! \!\!\!\!\! \!\!\bigoplus_{a\in\A,[w]_a \in W/R_a} \!\! \!\! \!\! \!\! \!\!
\str{M}\restr[w]_a\; , (P_a)_{a\in\A}\, ,\,
   R_I
    \bigr) 
  \]
  
We think of the disjoint union of the restrictions 
$\str{M}\restr[w]_a$
as realised over $a$-tagged versions of
$[w]_a$ with elements $u_a$ for $u \in [w]_a$ for the first sort
and $[w]_a$-tagged elements $s_\alpha$ for $s \in \wp([w]_a)$
in the second sort for $\alpha = [w]_a \in W/R_a$.  
More formally these can be 
realised over universes $[w_a] \times \{ a \}$ and
$\wp([w]_a) \times \{ \alpha \}$, but we use
notation like $u_a$ instead of $(u,a)$
for better readability.%
\footnote{One could also tag the elements $u_a = (u,a)$ for $u \in
  [w]_a$ slightly redundantly by their $\alpha = [w]_a \in W/R_a$;
  similarly,
  the $\alpha$-tag of $\emptyset \not= s \in \wp([w]_a)$ in
  $\str{M}\restr[w]_a $ is determined by $a$,
  since $u_a \in s$ implies $\alpha = [u]_a = [w]_a$,
  but we do distinguish $\emptyset \in \wp([w]_a)$ from
  $\emptyset \in \wp([w']_{a'})$ unless $a = a'$ and $[w']_{a'} = [w]_a
  \in W/R_a$.}
The interpretations for $E_a$ and the valuation $V$ are just restricted
over $[w]_a$ as indicated above.
The core structure $(W,(R_a))$ in the first-sort of
$\expl{\str{M}}$ is the underlying Kripke frame
on the set of worlds $W$ with the modal accessibility 
relations (representations of the $\sigma_a$ as
equivalence relations)
\[
R_a \;:=\; \{ (w,w') \in W^2 \colon  w' \in \sigma_a(w)   \}
\]
that make up this $S5$ frame; 
its domain $W$ also serves to supply the markers to
relate any $v \in W$ with its multiple copies in the local
structures  $\str{M}\restr [w]_a$ through $R_I$; 
the extra unary relation symbols $P_a$ carve out the subset of
$a$-tagged copies $w_a$ of worlds $w \in W$,
\[
\barr{r@{\;\;:=\;\;}l}
P_a &  \{ u_a \colon u \in [w]_a \},
\\
\hnt
R_I & \{ (w,w_a) \colon w \in W, a \in \A \}.
\earr
\]

One may read the index $I$ for ``incidence'' or
``identification''  of worlds with their multiple incarnations;
the subset $W$ of the first sort of $\expl{\str{M}}$ precisely serves
as a set of world labels that relate copies of the corresponding $w$:
\[
  u_a 
  \mbox{ represents the same world as }
  u'_{a'} 
  \mbox{ iff } (u_a, u'_{a'}) \in (R_I)^{-1} \!\!
  \circ R_I,
\]
i.e.\ iff there is some $v \in W$ for which $u_a, u'_{a'} \in
R_I[v]$, which implies that $v=u=u'$ and that
$[w]_a = [v]_a$ and $[w']_{a'} = [v]_{a'}$ for this $v = u = u'$.

The world-pointed model $\str{M},w$ is naturally represented by
the world-pointed exploded view, $\expl{\str{M}},w$ with distinguished
world~$w$ in the core of $\expl{\str{M}}$. In the 
state-pointed $\str{M},s$ with nonempty state $s \in S$, $s$ need not
be directly represented as an element of the second sort in any local
structure, as in general we only require $\wp(s) \subset S$ (cf.\
Definition~\ref{relencdef}).
But we may reduce this scenario to the world-pointed case by the
introduction of an extra (dummy) agent as follows.%
\footnote{The case of $\emptyset$-state-pointed
  models is trivial, since the empty information state is already trivially
  (and multiply) represented in the component structures $\str{M}\restr [w]_a$;
  so are all singleton information states.}
For non-trivial $s \subset W$ where $\wp(s) \subset S$ let
$\A_s := \A \,\dot{\cup}\, \{ s \}$ and consider the expansion of
$\str{M}$ to $\str{M}_s := (\str{M}, E_s)$ based on the relational
encoding of the trivial
inquisitive assignment
\[
  \Sigma_s \colon w \mapsto \left\{
    \barr{ll}
    \wp(s) & \mbox{ for } w \in s
    \\
    \wp(\{ w\}) & \mbox{ for } w \not\in s 
    \earr\right.
\]
whose induced $\sigma_s$ and $R_s$ partition $W$ into singleton
equivalence classes $[w]_s = \{ w \}$ for $w \notin s$ and the
single, potentially non-trivial equivalence class $[w]_s = s$, 
for any $w \in s$. So the exploded view of $\str{M}_s$ has, for non-trivial
$s$, an extra part for the $s$-local structure $\str{M}_s\restr [w]_s$
which provides a representation of $s$, both as a set of worlds in 
its first sort, and as the $\subseteq$-maximal element in its second sort.
So $s$ is $\FO$-definably available in $\expl{\str{M}_s}$  
in the $s$-local structure $\str{M}_s\restr [w]_s$, marked out either by 
the parameter $s = s_{\alpha_s}$, $\alpha_s := [w]_s$,
of the second sort, or by any choice of a parameter $w_s$, $w \in s$,
in the second sort.

We note that passage to $\str{M}_s$ is safe
also w.r.t.\ bisimulation 
equivalences for state-pointed models in the sense that
$\str{M},s \sim \str{M}',s'$ implies that $\str{M}_s,s \sim
\str{M}'_{s'},s'$, i.e.\ that for any $w \in s$ there is some $w' \in
s'$ such that $\str{M}_s,w \sim
\str{M}'_{s'},w'$ and vice versa; and similarly for every level of
$\sim^n$. So we give the following definition.

\bD
\label{explodedviewdef}
The \emph{exploded view} of a locally full relational epistemic model
$\str{M}$ is defined as $\expl{\str{M}}$ according to the above;
the analogous representation of an epistemic model $\M$ is obtained as
$\expl{\M} := \expl{\str{M}^\lf(\M)}$.
For world-pointed models we put
$\expl{\str{M},w} := \expl{\str{M}},w$.
In the non-trivially state-pointed version
we use $\expl{\str{M},s}:= \expl{\str{M}_s},s$ 
where $\str{M}_s$ is the expansion for a dummy agent $s$,
and $s$ 
stands for its representation $s_{\alpha_s}$ in the second sort of 
$\str{M}\restr \alpha_s$ as discussed above.
\eD

The exploded views $\expl{\str{M},w}$ and $\expl{\str{M},s}$ are such that the
evaluation of any $\FO$-formula in a single free first- or second-sort
variable over $\str{M},w$ or $\str{M},s$ translates into the
evaluation of a corresponding first-order formula over
$\expl{\str{M},w}$ and $\expl{\str{M},s}$, respectively. Indeed, there
is a uniform syntactic translation to this effect, a phenomenon 
that generally applies for first-order interpretations. In our case,
we rely on the uniform $\FO$ interpretability of $\str{M}$ within
its exploded view $\expl{\str{M}}$, together with the
$\FO$-definable access to the relevant assignment $w$ or $s$: 
$w \in W$ is an element of the
first sort of $\expl{\str{M},w}$; and $s \not=\emptyset$ is represented by the
maximal element in the second sort of its local 
component structure.%
\footnote{A standard technical detail involves the passage to quotients
w.r.t.\ $\FO$-definable equivalence relations for the representation
of especially the second sort of $\str{M}$ with $\expl{\str{M}}$.
By the very nature of an ``exploded view'', the
empty information state $\emptyset \in S$, all singleton states
$\{ w \} \in S$ and possibly other information states are represented in
the second sort of more than one local structure $\str{M}\restr[w]_a$.
The corresponding representation of the first sort $W$ is trivial
because $W$ is included as a definable subset in the first sort of
$\expl{\str{M}}$, even though structural properties of its elements
need to be retrieved from $R_I$-related copies in the
$\str{M}\restr[w]_a$ (see proofsketch below).}
In fact, for locally full inquisitive relational models $\str{M}$, 
the two-sorted structures $\str{M}$ and their exploded views
$\expl{\str{M}}$ are uniformly first-order interpretable
within one another (bi-interpretable). The following summarises
the use we want to make of interpretability in the context of
the upgrading scheme as indicated in Figure~\ref{S5upgradefigure}.

\bCl
\label{interpretationclaim}
Due to $\FO$-interpretability of any locally full relational model
$\str{M},w$ in $\expl{\str{M},w}$, or of 
$\str{M},s$ in $\expl{\str{M},s}$, there is a constant $d$ such that
$(q+d)$-elementary equivalence between the exploded views of two such
models implies $q$-elementary equivalence between
the original relational models. For instance in the world-pointed
case, any first-order \textup{$\phi(\w)$} of quantifier rank up to~$q$
in a world-variable has a uniform translation \textup{$\phi_I(\w)$}
of quantifier rank at most~$q+d$
such that $\str{M},w \models \phi$ iff
$\expl{\str{M},w} \models \phi_I$. It follows that  for
any pair of world-pointed locally full inquisitive relational models
$\str{M},w$ and $\str{M}',w'$
\[
  \expl{\str{M},w} \equiv_{q+d} \expl{\str{M}',w'} \quad \Rightarrow
\quad
\str{M},w \equiv_q \str{M}',w';
\]
and similarly for formulae \textup{$\phi(\mathtt{s})$} in a state variable and
state-pointed locally full relational epistemic models. 
\eCl

\prf[Proofsketch.]
All the relevant structural information about $\str{M}$ is easily imported to
the sub-universe $W$ of the first sort in $\expl{\str{M}}$. This
sub-universe itself is $\FO$ definable in $\expl{\str{M}}$
as the projection of $R_I$ to its first component; the valuation $V$
at $w \in W$ is imported from any one of the elements in $R_I[w]$;  
the relational encoding $E_a[w]$ of the inquisitive assignment
$\Sigma_a(w)$ is similarly imported through backward translation of $E_a(w_a)$
from any one of the $w_a \in R_I[w] \cap P_a$ via set-wise passage to
$R_I$-predecessors:
\[
  v \in s \in \Sigma_a(w) \quad\mbox{ iff }\quad
  v \in (R_I \circ \in)^{-1}(E_a[w_a])
\]
for $w_a \in R_I[w] \cap P_a$.
Compare Figure~\ref{explfigure} to
see how especially this relationship between elements and subsets of the core
part can be traced in the representation of the $a$-local structure
$\str{M}\restr[w]_a$ in the halo of $\expl{\str{M}}$. 
All these relations are therefore uniformly first-order definable
within the exploded views $\expl{\str{M}}$, and the maximal 
quantifier rank in corresponding $\FO$-formulae can serve as the value
$d$ for the transfer claim.
For the non-trivially state-pointed case, similarly, the distinguished state $s$ (and
its entire powerset) is accessible in the component structure
$\str{M}_s\restr [w]_s$ for any $w \in s$ in $\expl{\str{M},s}$. 
\eprf

For the upgrading argument according to
Figure~\ref{S5upgradefigure}
we shall therefore eventually establish
$\equiv_q$-equivalence of, for instance world-pointed inquisitive relational models,
\begin{equation}
   \str{M},w \equiv_q \str{M}',w'
  \label{upgrade1}
\end{equation}
where $q$ is the quantifier rank of the given bisimulation invariant
$\phi(x) \in \FO$, based on the assumption that
$\str{M},w \sim^\ell \str{M}',w'$
for suitable $\ell =
\ell(q)$. The above preparation allows us to reduce this to showing
\begin{equation}
  \expl{\str{M},w}= \expl{\str{M}},w \equiv_{q+d} \expl{\str{M}'},w'
  = \expl{\str{M}',w'}.
  \label{upgrade2}
\end{equation}

In fact we also only need to establish~(\ref{upgrade1}) or~(\ref{upgrade2})
for $\str{M}$ and $\str{M}'$
from some subclass of the class of all 
epistemic relational models that represents
the corresponding class (of all or just all finite such models)
up to bisimulation equivalence. E.g.\ in the finite model theory
version for the world-pointed case, it suffices to show that, for any two finite
epistemic models
$\M,w$ and $\M',w'$ that are $\ell$-bisimilar,
there are bisimilar companions
$\hat{\M},\hat{w} \sim \M,w$
and $\hat{\M}',\hat{w}' \sim \M',w'$
such that 
\begin{equation}
  \expl{\hat{\M},\hat{w}}
= \expl{\str{M}^\lf(\hat{\M})},\hat{w}
\equiv_{q+d}
\expl{\str{M}^\lf(\hat{\M}')},\hat{w}' =
\expl{\hat{\M}',\hat{w}'}
\label{upgrade3}
\end{equation}

The remaining work towards this upgrading argument, therefore, focuses on
the construction of suitable bisimilar companions. Suitability here
means that (\ref{upgrade3}) can be established by 
Ehrenfeucht--Fra\"\i ss\'e arguments
for $\FO$ and -- for auxiliary arguments concerning the
restrictions $\str{M} \restr [w]_a$ as main building blocks of the
$\expl{\str{M}}$ --  for $\MSO$.

The Ehrenfeucht--Fra\"\i ss\'e argument  for suitably prepared 
epistemic models is facilitated by the observation that 
$\sim$-invariance implies invariance under
disjoint unions for epistemic models,
and that the passage to exploded views is compatible with such disjoint unions.
The following definition of disjoint unions is the natural one for
inquisitive models.

\bD
\label{disjointuniondef}
The \emph{disjoint union} of epistemic models $\M_1$, $\M_2$, denoted as
$\M_1 \oplus \M_2$, is the epistemic model whose set of worlds is the
disjoint union of the sets of worlds $W_i$ of the $\M_i$, with the
obvious transfer of inquisitive and propositional assignments from
the $\M_i$ to   $\M_1 \oplus \M_2$. 
\eD

Clearly, $\M,w \sim \M \,\oplus\, \mathbb{L},w$ and
$\M,s \sim \M \,\oplus\, \mathbb{L},s$ for any world- or state-pointed epistemic
models $\M,w$ or $\M,s$ and plain $\mathbb{L}$. But while the disjoint union of
epistemic models or of their relational encodings cannot be strictly
disjoint, because of the special r\^ole of the empty information state
$\emptyset$, this little problem is overcome in their exploded views.

\bO
\label{disjunionobs}
Passage to the exploded view of epistemic models
is compatible with disjoint unions in the sense that
\[
 \expl{\M_1 \oplus \M_2}=
 \expl{\M_1}  
 \oplus \expl{\M_2},
\]
and similarly with a distinguished world or state parameter from
$\M_1$ say.
\eO

\paragraph*{\boldmath FO-equivalence and Gaifman locality.} 
We review some relevant technical background on Gaifman locality 
and its use towards establishing levels of first-order
equivalence between relational structures, with a view to applying these techniques
to show $\expl{\hat{\str{M}}, \cdot\, } \equiv_{q+d} \expl{\hat{\str{M}}' , \cdot\, }$
in the situation of Figure~\ref{S5upgradefigure}, for suitably
pre-processed models $\hat{\str{M}}, \cdot\, \,\simn\, \hat{\str{M}}' , \cdot\, $.
The relational structures in question interpret just unary and binary 
relations. Correspondingly, \emph{Gaifman distance} over their (two-sorted)
universes is just ordinary graph distance w.r.t.\ the symmetrisations
of the binary relations.
This establishes a natural distance measure between elements, in our case across
both sorts, worlds and information states.%
\footnote{In general, `infinite distances' occur between elements
  that are not linked by finite chains of binary relational edges,
  i.e.\ that are in distinct connected components.}
For an element $b$ of the relational structure $\str{B}$ and $\ell \in \N$, 
the set of elements at distance up to $\ell$ from $b$ is the
\emph{$\ell$-neighbourhood} of $b$, denoted as 
\[
N^\ell(b) = \{ b' \in B \colon d(b,b') \leq \ell \}.
\]  

A first-order formula $\psi(x)$ in a
single free variable $x$ (of either sort) is \emph{$\ell$-local} if its semantics at $b \in
\str{B}$ depends just on the induced substructure 
$\str{B}\restr N^\ell(b)$ on the $\ell$-neighbourhood $N^\ell(b)$.
In other words, $\psi(x)$ is $\ell$-local if, and only if, it is
logically equivalent to its 
\emph{relativisation} to the $\ell$-neighbourhood  of $x$,
i.e., to the formula $\psi^\ell(x)$ obtained by explicitly restricting all quantifiers
in $\psi(x)$ to the $\ell$-neighbourhood of $x$ (which is possible provided the
underlying relational signature is finite).
We refer to $\psi^\ell(x) \in \FO$ as the 
\emph{$\ell$-localisation} of $\psi(x)$.
A finite set of elements of the relational
structure $\str{B}$ is said to be \emph{$\ell$-scattered} if the
$\ell$-neighbourhoods of any two distinct members are disjoint. 
A \emph{basic local sentence} is a sentence that asserts, for some
formula $\psi(x)$ and some $m \geq 1$, 
the existence of an $\ell$-scattered set of $m$ elements that 
each satisfy the $\ell$-localisation $\psi^\ell(x)$. 
Gaifman's theorem~\cite{Gaifman,EF}, in the case that we are
interested in, asserts that any first-order formula in a single free
variable $x$ is logically equivalent to a boolean combination 
of $\ell$-local formulae $\psi^\ell(x)$ for some $\ell$ and some basic
local sentences. We shall focus on the following
levels of Gaifman equivalence.

\bD
\label{Gaifmandef}
Two pointed relational structures $\str{B},b$ and $\str{B}',b'$ are
\emph{$\ell$-locally $r$-equivalent}, denoted
$\str{B},b \equiv_r^{\ssc (\ell)} \str{B}',b'$, if $b$ and $b'$ satisfy
exactly the same $\ell$-localisations of formulae $\psi(x)$ 
of quantifier rank up to~$r$.

$\str{B},b$ and $\str{B}',b'$ are 
\emph{$(\ell,r,m)$-Gaifman-equivalent}, denoted
$\str{B},b \equiv_{r,m}^{\ssc (\ell)} \str{B}',b'$, if
$\str{B},b$ and $\str{B}',b'$ are
$\ell$-locally $r$-equivalent and
$\str{B}$ and $\str{B}'$ satisfy exactly the same 
basic local sentences concerning $\ell'$-localisations  
for formulae $\psi(x)$ of quantifier rank up to~$r$ and $\ell'$-scattered
sets of sizes up to~$m$ for any $\ell' \leq \ell$.
\eD

The following is an immediate corollary of 
Gaifman's classical theorem~\cite{Gaifman,EF}, also in the case of
our two-sorted structures and for distinguished elements of either
sort. 

\begin{figure}
  \[
    \nt\hspace{-1cm}
\xymatrix{
\str{M},w \ar@{-}[d]|{\rule{0pt}{1ex}\sim}\ar @{-}[rr]|{\;\simn\,}
&& \str{M}'\!,w' \ar @{-}[d]|{\rule{0pt}{1ex}\sim}
\\
\hat{\str{M}},\hat{w} \ar@{-}[d]|{\rule{0pt}{1ex}\sim}\ar @{-}[rr]|{\;\simn\,}
&& \hat{\str{M}}'\!,\hat{w}' \ar @{-}[d]|{\rule{0pt}{1ex}\sim}
\\
\hat{\str{M}},\hat{w} \, \oplus \str{X} 
\ar@{-}[rr]|{\rule{0pt}{1.5ex}\;\equiv_q}
\ar@{<->}[d]_{\expl{\cdot}}
&& \hat{\str{M}}',\hat{w}' \,\oplus \str{X}  
\ar@{<->}[d]^{\expl{\cdot}}
\\
\expl{\hat{\str{M}},\hat{w}
\,\oplus \str{X}}  
\ar@{-}[rr]|{\rule{0pt}{1.5ex}\; \equiv^{\ssc (\ell)}_r}
&& \expl{\hat{\str{M}}',\hat{w}' \,\oplus \str{X}}  
\\
\expl{\hat{\str{M}},\hat{w}}  
\ar@{-}[rr]|{\rule{0pt}{1.5ex}\; \equiv^{\ssc (\ell)}_r}
\ar@{->}[u]|{\rotatebox{90}{$\ssc \,\subset\,$}}
&& \expl{\hat{\str{M}}',\hat{w}'} 
\ar@{->}[u]|{\rotatebox{90}{$\ssc \,\subset\,$}}
}
\qquad\qquad
\xymatrix{
\str{M}_s,s \ar@{-}[d]|{\rule{0pt}{1ex}\sim} \ar @{-}[rr]|{\;\simn\,}
&& \str{M}'_{s'}\!,s' \ar @{-}[d]|{\rule{0pt}{1ex}\sim}
\\
\hat{\str{M}}
,\hat{s} \ar@{-}[d]|{\rule{0pt}{1ex}\sim}\ar @{-}[rr]|{\;\simn\,}
&& \hat{\str{M}}'
\!,\hat{s}' \ar @{-}[d]|{\rule{0pt}{1ex}\sim}
\\
\hat{\str{M}}
,\hat{s} \,\oplus \str{X}  
\ar@{-}[rr]|{\rule{0pt}{1.5ex}\;\equiv_q}
\ar@{<->}[d]_{\expl{\cdot}}
&& \hat{\str{M}}'
\!,\hat{s}' \,\oplus \str{X}  
\ar@{<->}[d]^{\expl{\cdot}}
\\
\expl{\hat{\str{M}}
  ,\hat{s} \,\oplus \str{X}} 
\ar@{-}[rr]|{\rule{0pt}{1.5ex}\; \equiv^{\ssc (\ell)}_r}
&& \expl{\hat{\str{M}}'
  \!,\hat{s}' \,\oplus \str{X}} 
\\
\expl{\hat{\str{M}},\hat{s}}  
\ar@{-}[rr]|{\rule{0pt}{1.5ex}\; \equiv^{\ssc (\ell)}_r}
\ar@{->}[u]|{\rotatebox{90}{$\ssc \,\subset\,$}}
&& \expl{\hat{\str{M}}',\hat{s}'}  
\ar@{->}[u]|{\rotatebox{90}{$\ssc \,\subset\,$}}
}
\]
\caption{Upgrading for relational epistemic models,
  refined by Proposition~\protect\ref{Gaifmanprop}.}
\label{S5extendedupgradefigure}
\end{figure}
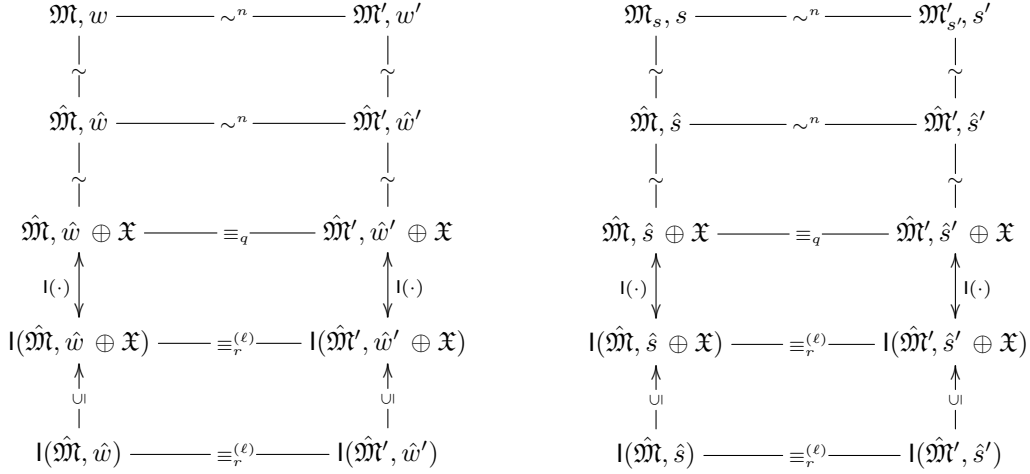

\bP
\label{Gaifmanprop}
The semantics of any first-order formula $\psi(x)$ in a finite purely relational vocabulary  
is preserved under Gaifman equivalence $\equiv_{r,m}^{\ssc (\ell)}$
for sufficiently large values of the parameters $\ell, r,m \in \N$.
For formulae $\psi(x) \in \FO$ whose semantics is also preserved under
disjoint unions of relational structures, preservation under 
$\equiv_{r,m}^{\ssc (\ell)}$ implies preservation under
$\equiv_r^{\ssc (\ell)}$.
\eP

We can use these results from classical first-order model theory in
order to simplify the upgrading task: in effect we may replace the
target levels $\equiv_q$ from Figure~\ref{genupgradefigure}, or 
$\equiv_{q+d}$ in Figure~\ref{S5upgradefigure}
by $\equiv_r^{\ssc (\ell)}$ for exploded views,
for appropriate $\ell$ and $r$. 
Compare Figure~\ref{S5extendedupgradefigure} for the following argument. 
On the basis of interpretability, $\equiv_q$ over
$\hat{\str{M}}/\hat{\str{M}}'$ follows from
$\equiv_{q+d}$ over $\expl{\hat{\str{M}}}/\expl{\hat{\str{M}}}'$
as illustrated in Figure~\ref{S5upgradefigure}, which in turn is
implied by a suitable level of $\equiv_{r,m}^{\ssc (\ell)}$ over these.
Up to bisimilarity, we may replace models
$\hat{\str{M}}$ and $\hat{\str{M}}'$ by disjoint sums with a model
\[
\str{X} \;:=\; m \!\otimes\! \hat{\str{M}} \;\oplus\; m\!\otimes\!\hat{\str{M}}'
\]
consisting of $m$ many disjoint copies of each of $\hat{\str{M}}$ and
$\hat{\str{M}}'$ (and similarly, using $\hat{\str{M}}_{\hat{s}}$ and
$\hat{\str{M}}'_{\hat{s}'}$ in the non-trivially state-pointed case).
The purpose of this modification is simply to
trivialise the contribution of basic local sentences about scattered
sets of size up to~$m$ towards $\equiv_{r,m}^{\ssc (\ell)}$.   
By compatibility of the exploded view with disjoint
sums,
$\equiv_{r,m}^{\ssc (\ell)}$ and $\equiv_{r}^{\ssc (\ell)}$ coincide
for the exploded views of these augmented structures,
which further reduces the issue to the $\ell$-local equivalence  
\begin{equation}
  \expl{\hat{\str{M}},\cdot}
\;\equiv_{r}^{\ssc (\ell)}\;
\expl{\hat{\str{M}}',\cdot}
\label{Gaifmanequiveqn}
\end{equation}
in both the world-pointed and the state-pointed case, as illustrated
in Figure~\ref{S5extendedupgradefigure}.

The following summarises these
considerations, which pave
the way towards establishing 
the compactness property which we have seen
to be equivalent to our main theorem.
We can always assume that the underlying signature
is the signature of the given first-order formula, hence finite. 

\bP
\label{summaryGaifmanupgradeprop}
In order to show that any first-order formula in a single free world or
state variable that is $\sim$-invariant over a class $\CC$ of world-
or state-pointed relational epistemic models is even invariant under 
$\simn$ for some finite level~$n \in \N$, and hence expressible in
$\inqbm$ over $\CC$, it suffices to show that, for any $\ell, r \in \N$
there is an $n \in \N$ such that any two $\simn$-equivalent models
from $\CC$ can, up to $\sim$, be replaced within $\CC$ by models,
whose exploded views are $\ell$-locally $\FO$-equivalent
up to quantifier rank~$r$.
\eP

Towards the criteria in Proposition~\ref{summaryGaifmanupgradeprop},
we exploit the representation of the relevant models in their
exploded views. These allow us to view an inquisitive
epistemic model $\M$ as a network of interconnected $a$-local structures
$\M\restr [w]_a$ that are linked through shared worlds.
The individual mono-modal $a$-local structures, where two-sortedness reflects the
higher-order nature of the inquisitive assignment, can be managed
well in isolation. The interaction between these $a$-local structures,
for different agents~$a$ in the global pattern of overlapping
$a$-classes, on the other hand, is essentially reflected in
the single-sorted $S5$-Kripke structure $\str{K}(\M)$ that forms
the core part of the exploded view $\expl{\M}$ in the  first sort.
The desired upgrading argument for
Proposition~\ref{summaryGaifmanupgradeprop}
calls for a strategy in the $r$-round first-order
Ehrenfeucht--Fra\"\i ss\'e\ game over the $\ell$-neighbourhoods
of exploded views of two suitable pre-processed $n$-bisimilar models.
Their nature as conglomerates of interconnected $a$-local structures
allows us to divide concerns strategically as follows: 
\begin{itemize}
\item[--]
  $a$-local concerns in suitably regularised local structures:
  the analysis of the $\FO$-game over these two-sorted structures
  largely reduces to $\MSO$-analysis over their first sorts;
  local pre-processing towards a normalisation of the inquisitive
  assignments will support that analysis.
  \item[--] the actual Ehrenfeucht--Fra\"\i ss\'e game in the $\ell$-neighbourhoods, 
involving the connectivity pattern of links between $a$-local component
structures; for that, global pre-processing towards local acyclicity
in the connectivity patterns supports an analysis based on tree-like hypergraphs.  
\end{itemize}
And indeed, the structure of our exploded views is precisely designed for this
purpose. Due to suitable global pre-processing the core parts will be
locally tree-like, and the general format of the exploded view as
illustrated in Figure~\ref{explfigure} retains this tree-likeness 
because local structures are connected only through the core; 
and the tentacles into the second sort, which each belong to just one
local structure, will have been tamed by local pre-processing.
Both the local and the global pre-processing need to respect
bisimulation equivalence, i.e.\ rely on the construction of suitable bisimilar companions
within the class~$\CC$ at hand.

\subsection{Suitable bisimilar companions}
\label{FOscattersec}

To deal with various uniform constructions of bisimilar companions
we use the notion of (globally bisimilar) coverings as expounded e.g.\ 
in~\cite{OttoAPAL04,DawarOttoAPAL09,Otto12JACM}.
The natural adaptation to the setting of multi-modal inquisitive 
epistemic models is the following.

\bD
\label{inqbisimcoverdef}
A \emph{bisimilar covering} of an inquisitive model
$\M = (W, (\Sigma_a),V)$ 
by an inquisitive model
$\hat{\M} = (\hat{W}, (\hat{\Sigma}_a),\hat{V})$
is a map
\[
\pi \colon \hat{\M} \longrightarrow \M,
\]
based on a surjection $\pi \colon \hat{W} \rightarrow W$ 
with natural induced maps $\pi \colon \wp(\hat{W})
\rightarrow \wp(W)$ and $\pi \colon \wp\wp(\hat{W}) \rightarrow
\wp\wp(W)$, 
such that 
$\pi$ is compatible with the inquisitive and propositional assignments
in the sense that the following diagram commutes
\[
\xymatrix{
\hat{W} \ar[d]^\pi  \ar[rr]^{\hat{\Sigma}_a} 
&& \wp\wp(\hat{W}) \ar[d]^\pi 
\\
W  \ar[rr]^{\Sigma_a} && 
\wp\wp(W)
}
\]
and that $\hat{V}(p) = \pi^{-1}(V(p))$ for all $p \in \P$
(i.e.\ $\pi$ is a strict homomorphism w.r.t.\ the valuations).
The covering is \emph{finite} if all its fibres
$\pi^{-1}(w) \subset \hat{W}$ are finite. 
\eD

The definition implies that $\pi$ induces a global bisimulation
between $\M$ and $\hat{\M}$
as $\hat{\M},\hat{w} \sim \M,\pi(\hat{w})$
for all $\hat{w} \in \hat{W}$.

Towards the desired upgrading arguments,
a partial unfolding of the global pattern of 
equivalence classes $[w]_a$ in an inquisitive $S5$ model can be
achieved in finite bisimilar coverings.
The eventual covering construction at the global level
requires a pre-processing of the local component
structures $\str{M}\restr [w]_a$ of $\expl{\str{M}}$ as introduced
in Section~\ref{FOEFGaifsec} above:
\[
\str{M}\restr[w]_a:= \bigl( [w]_a, \wp([w]_a), E_a \restr[w]_a, (P_i \restr[w]_a) \bigr),
\]
which for a locally full $\str{M}$ is the full relational encoding
of the mono-modal, single-component epistemic model
\[
  \M\restr[w]_a=
  ([w]_a, \Sigma_a\restr [w]_a, V \restr [w]_a),
\]
whose epistemic
inquisitive assignment $\Sigma_a \colon [w]_a \rightarrow \wp\wp([w]_a)$
is constant and can be identified with its value
$\Sigma_a(w) \in \wp\wp([w]_a)$.

\medskip
We first explore the extent to which this $\Sigma_a$-assignment can be
locally modified and simplified in companion models that retain
bisimilarity in the global context. 
This  local pre-processing in Section~\ref{localpreprocsec} is specific for
epistemic $\inqbm$ and new here; the techniques for the global pre-processing
via local unfoldings in globally bisimilar finite covers, in Section~\ref{epicovsec},
adapts ideas from modal logics in~\cite{OttoAPAL04,DawarOttoAPAL09}; 
both aspects go substantially beyond the treatment of basic $\inqbm$ in~\cite{CiardelliOttoJSL2021}.

\subsubsection{Local pre-processing.}
\label{localpreprocsec}
Since $\Sigma_a$ is constant across the single equivalence class $[w]_a$,
bisimulation
-- for this single-agent structure --
equivalence collapses to its second level $\sim^1$.
W.r.t.\ a fixed set $\P$ of basic propositions 
each bisimulation type of such a structure is characterised
by the collection of propositional types realised overall in its
worlds and in the information states in its inquisitive 
assignment. In particular, multiplicities of these propositional
types do not matter at all, neither overall nor in any individual
information state: any positive number of realisations is as good as
just a single one, so that multiplicities can be uniformly boosted to
avoid distinguishing features at the level of $\FO$.

\medskip
The natural restrictions-cum-reducts
$\M\restr[w]_a$ ignore the inquisitive structure
induced by agents $b \not= a$, and fall into one of these 
global
$\sim^1$-equivalence classes.
All the remaining complexity of multi-agent inquisitive $S5$
structures arises from the overlapping of $a$-classes
for different $a \in \A$, i.e.\ from the manner in which the various
local structures $\M\restr [u]_a$ or their relational encodings
$\str{M}\restr [w]_a$ overlap. 
In order to freeze the manner in which the local
structures on individual $a$-classes are stitched together, we 
endow the $\M\restr [w]_a$ with colourings of their worlds that
prescribe their bisimulation types in the global structure $\M$.
For different purposes we use different
granularities for this colouring. The finest one of interest 
is the colouring based on the full bisimulation-type in $\M$, with colour set 
$C = W/{\sim}$, the set of $\sim$-equivalence classes of worlds in
$\M$. This is the colouring that needs to be respected in passage to fully bisimilar
  partner structures like $\hat{\str{M}},\hat{w} \sim \str{M},w$ in
  the vertical axes in Figure~\ref{S5extendedupgradefigure}.

Coarser colourings are induced by the $\simm$-type, for a fixed finite
level $m$, with colour set $C = W/{\simm}$. These may also be 
obtained via the natural projection that identifies all
$\sim$-types that fall into the same $\simm$-class.
Successive coarsenings of this kind will inform the
strategy analysis towards establishing levels of
 $\FO$-indistinguishability like
$\hat{\str{M}},\hat{w}
\,\oplus \str{X} \; \equiv^{\ssc (\ell)}_r \hat{\str{M}}',\hat{w}'
\,\oplus \str{X}$ on the basis of $\str{M},w \sim^n \str{M}',w'$
across the horizontal axes in Figure~\ref{S5extendedupgradefigure}.
Moreover, the relevant granularity~$m$ will be dependent on the
distance from the distinguished worlds or states, just as control over
levels of bisimilarity decreases with the depth of the exploration.
But this granularity~$m$ will be uniformly bounded by the initial
level~$n$ of $n$-bisimilarity in $\str{M},w \sim^n \str{M}',w'$ or $\str{M},s \sim^n
\str{M}',s'$. 

Formally, a colouring of worlds by  any colour set $C$ can be encoded
as a propositional assignment, for new atomic 
propositions $c \in C$. Unlike
propositional assignments in general, these propositions will be mutually
exclusive, so as to partition the set of worlds.
To reduce conceptual overhead, we
encode colourings as functions
\[
\barr[t]{rcl}
\rho \colon W &\rightarrow& C 
\\
w &\mapsto& \rho(w) 
\earr
\quad \mbox{ instead of } \quad
\barr[t]{rcl}
V \colon C &\rightarrow& \wp(W)
\\
c &\mapsto& V(c) = \{w \colon \rho(w) = c\}
\earr
\]
and correspondingly write, e.g., 
$(W, (\Sigma_a)_{a \in \A},
\rho)$
for the inquisitive epistemic model with 
the colouring $\rho \colon W \rightarrow C$ that assigns
to each world its $\sim$-class: 
\[
\barr[t]{rcl}
\rho \colon \; W &\longrightarrow& C := W/{\sim} 
\\
\hnt
w & \longmapsto & [w]_{\sim} 
= \{ u \in W \colon \M,u \sim \M,w \}
\earr 
\]

The coarser colouring based on $\simm$-types 
is correspondingly formalised in 
the inquisitive epistemic model $(\M, (\Sigma_a)_{a \in \A},
\rho_m)$, with 
the colouring $\rho_m \colon W \rightarrow W/{\simm}$.
Clearly the $\sim$-colouring $\rho$ refines the $\simm$-colouring $\rho_m$, and 
$\rho_n$ refines $\rho_m$ for $m<n$. 
In particular, all levels refine $\rho_0$, which determines each 
world's original propositional type as induced by $V$ in $\M$.

Finiteness of the set of colours $W/{\sim}$ 
for the $\sim$-colouring $\rho$ is 
guaranteed over each individual finite model $\M$ (or across pairs of finite 
models $\M$ and $\M'$ to be compared). But also over possibly infinite 
models the sets of colurs  $W/{\simm}$ for the $\simm$-colourings $\rho_m$
are unifoirmly finitely bounded across
all models in a finite signature, for each finite level~$m$.

\bD
\label{localastrgranulardef}
For an $a$-class $[w]_a$ in an inquisitive model $\M$,
let the associated \emph{local $\sim$-structure} 
be the mono-modal inquisitive $S5$ structure $\M\restr [w]_a,\rho$ 
with the propositional assignment induced by $\rho$ in  $\M$, which
colours each world in $[w]_a$ according to its $\sim$-type in $\M$.

For $m \in \N$, let $\M\restr [w]_a,\rho_m$ be the coarsening
of granularity $m$ of that local structure:
the \emph{local $\simm$-structure} with 
propositional assignment induced by $\rho_m$, which
colours each world in $[w]_a$ by its $\simm$-type in $\M$
according to the natural projection $\pi \colon W/{\sim} \rightarrow W/{\sim^m}$.
\eD

Interestingly, a local $\sim$- or $\simm$-structure can be modified
considerably without changing its bisimulation type.
We now want to achieve local structures 
\begin{itemize}
\item[--]
that realise every bisimulation type of worlds locally with high
multiplicity; this leads to the notion of \emph{richness} in 
Definition~\ref{richdef};
\item[--]
  whose inquisitive assignment $\Sigma_a$ is normalised w.r.t.\ 
  to the structure
  of the inquisitive state $\Pi = \Sigma_a(w)$ within $\wp([w]_a)$ so
  as to suppress $\FO$-definable distinctions that are not governed
  by any level of $\simm$;
  this leads to the notion of
  \emph{regularity} in Definition~\ref{regulardef}.
\end{itemize}

\paragraph*{Richness.} 
The idea for \emph{richness} is a
quantitative one, viz., that
all bisimulation types that are realised in any information
state in $\Sigma_a(w)$ 
will be realised with high multiplicity in some superset that is
also in $\Sigma_a(w)$. Any desired finite level of richness
can be achieved in finite bisimilar coverings
(Lemma~\ref{richcoverlem}), which 
essentially just put a fixed number of copies of every world.

\bD
\label{richdef}
A local structure $\M \restr [w]_a$ 
in $\M$ is \emph{$K$-rich} for some $K \in \N$ if,
for every information state $s \in \Sigma_a(w)$ there is some 
$s'$, $s \subset s' \in \Sigma_a(w)$, such that
every $\sim$-type in $W/{\sim}$ that is realised in
$s$ is realised at least $K$~times in~$s'$.
An inquisitive $S5$ structure $\M = (W,(\Sigma_a),V)$ is 
\emph{$K$-rich} if each one of its local 
structures is $K$-rich, for all $a \in \A$.
\eD

\bL
\label{richcoverlem} 
For any $K \in \N$, any inquisitive $S5$ structure $\M = (W, (\Sigma_a),V)$
admits a finite bisimilar covering that is $K$-rich. If finiteness
does not matter, richness can similarly be achieved for any target
cardinality.  
\eL

\prf
It suffices to take the natural product $\M \times [K]$ of $\M$ with the standard 
$K$-element set $[K] = \{1,\ldots, K \}$. This results in a model whose
 set of worlds is
$W \times [K]$, whose propositional assignment is induced by the projection 
$\pi \colon W \times [K] \rightarrow W$, and whose inquisitive
assignment  
is $\Sigma_a(w,m) := 
\{ s \subset W \times [K] \colon \pi(s) \in \Sigma_a(w) \}$.
It is easily checked that $\M \times [K]$ is $K$-rich and that $\pi$ is a
bisimilar covering in the sense of Definition~\ref{inqbisimcoverdef}.
\eprf

\paragraph{Regularity.}
Unlike the merely quantitative flavour of richness, 
regularity imposes qualitative constraints on the structure of
the family of information states in $\Sigma_a(w)$. It requires a
modification of the local inquisitive assignment and 
is not obtained in a covering. The idea is 
to normalise the downward closed family of information states 
in $\Pi := \Sigma_a(w)$ within the boolean algebra $\wp([w]_a)$ of subsets of
$[w]_a$ as much as possible. This normalisation needs to preserve 
the bisimulation type as encoded in the colouring $\rho$.  
It must also be geared towards passage to the coarser picture 
induced by increasingly lower levels of $m$-bisimulation as encoded in
the $\rho_m$.%
\footnote{An oversimplification in earlier
  versions of this preprint wrongly assumed robustness w.r.t.\ such
  coarsenings; we are grateful to an anonymous referee for spotting
  the mistake.}
Passage to some suitable degree of regularity  
should be thought of as a local modification of the representation,
rather than the epistemic content, of the inquisitive assignment. 
Indeed, bisimilarity is rather 
robust under changes in the actual composition of the inquisitive
assignment $\Pi = \Sigma_a(w)  \subset \wp([w]_a)$ in $\M$.
By the inquisitive epistemic nature of
$\M$, $\Sigma_a$ has constant value $\Pi$ on $[w]_a$
and $\Pi$ is a downward closed collection of subsets of $[w]_a$
such that $w \in \bigcup \Pi =  [w]_a = \sigma_a(w)$. The only additional invariant
imposed on $\Pi$ up to $\sim$ stems from the
associated colouring $\rho$ of worlds by their $\sim$-types
in $C = W/{\sim}$; for the inquisitive assignment $\Pi$ this invariant
concerns the colour combinations in 
\[
\barr{r@{}l}
\rho(\Pi)  & \;:=\; \{ \rho(s) \colon s \in \Pi \} 
\subset \wp(C) 
\\
\hnt
\mbox{ for }\; &
\rho(s) := 
\{ \rho(v) \colon v \in s \} = 
\{ c \in W/{\sim} \colon c \cap s
\not= \emptyset \}.
\earr
\]
 
Note that $\rho(\Pi)$ is downward closed in $\wp(C)$.
Not only the collection of bisimulation types in the set of worlds 
$[w]_a$, but also in the information states of the inquisitive assignment $\Pi$
is fully determined by $\rho(\Pi)$.
Up to bisimulation, $\Pi$ can be replaced by any inquisitive 
epistemic assignment $\hat{\Pi}$ with $\rho(\hat{\Pi}) = \rho(\Pi)$.

Similar remarks apply to finite levels $\simm$ relative to the
colouring $\rho_m$ with values in $W/{\simm}$,
and for analogously defined $\rho_m(\Pi)= \{
\rho_m(s) \colon s \in \Pi \}$. Clearly $\rho_m(\Pi)$ is fully
determined by $\rho(\Pi)$, so that $\rho_m(\hat{\Pi}) = \rho_m(\Pi)$
for all admissible variations $\hat{\Pi}$.
But crucially, due to the coarsening,
colour combinations $\rho_m(s_i)$ for infomation states
$s_1,s_2 \in \Pi$ may coincide even if $\rho(s_1) \not= \rho(s_2)$. 
So the constraint that $\rho(\hat{\Pi}) = \rho(\Pi)$ implies that some
colour combinations in $\rho_m(\hat{\Pi}) = \rho_m(\Pi)$ may have to
appear in more than one $\subset$-maximal element of $\rho_m(\Pi)$
and hence must be realised in at least that many $\subset$-maximal
elements of any admissible $\hat{\Pi}$.%
\addtocounter{footnote}{-1}\footnotemark

The assignment $\hat{\Pi}$ we choose
will therefore be obtained from an upper bound $\kappa$ for the possibly unavoidable
multiplicity of $\subset$-maximal elements in $\rho_m(\Pi)$.
So $\kappa$ is to be adjusted accross the two models to be matched
in the upgrading argument of Figure~\ref{S5extendedupgradefigure}; 
it will always be finite in the context of 
finite inquisitive models for the finite model theoretic reading of
our results, but could necessarily be infinite in the general case of infinite models that
might even locally realise infinitely many distinct bisimulation types. 
Its granularity levels~$m$ will be dependent on the distance of
the local structures from the distinguished worlds or states in the upgrading
scenario of Figure~\ref{S5extendedupgradefigure}. It will always be
uniformly bounded by the initial level of $n$-bisimilarity $\sim^n$
between the two pointed
models under consideration. So there is a fixed finite bound on the
number of $\simm$-colours that may occur. 
Such an upper bound can uniformly be chosen as a bound on
the number of distinct $\sim$-types realised by information
states in the inquisitive assignments under consideration.

\begin{figure}
$\color{blue}
\underbrace{%
{\begin{tikzpicture}[scale=.8]
\node at (-2.8,-2.5) [color=blue, opacity=.6]  {$B_1$};
\node at (-2.6,-1) [color=blue, opacity=.7]{\rotatebox{90}{$s(\alpha_1,1)$}};
\node at (-1.5,-1) [color=blue, opacity=.7]{\rotatebox{90}{$s(\alpha_2,1)$}};
\node at (.7,-1) [color=blue, opacity=.7]{\rotatebox{90}{$s(\alpha_n,1)$}};
\node at (-.4,-1) [color=blue, opacity=.4] {{\boldmath \large $\cdots$}};
\path[fill=blue, fill opacity=.1] (-3.1,-2) rectangle ++(1cm,2cm);
\path[fill=blue, fill opacity=.1] (-2,-2) rectangle ++(1cm,2cm);
\path[fill=blue, fill opacity=.1] (.2,-2) rectangle ++(1cm,2cm);
\path[fill=blue, fill opacity=.1] (-3.2,-2.1) rectangle
++(4.5cm,2.2cm);
\end{tikzpicture}}
\;
{\begin{tikzpicture}[scale=.8]
\node at (-2.8,-2.5) [color=blue, opacity=.6]  {$B_2$};
\node at (-2.6,-1) [color=blue, opacity=.7]{\rotatebox{90}{$s(\alpha_1,2)$}};
\node at (-1.5,-1) [color=blue,opacity=.7]{\rotatebox{90}{$s(\alpha_2,2)$}};
\node at (.7,-1) [color=blue, opacity=.7]{\rotatebox{90}{$s(\alpha_n,2)$}};
\node at (-.4,-1) [color=blue, opacity=.4] {{\boldmath \large $\cdots$}};
\path[fill=blue, fill opacity=.1] (-3.1,-2) rectangle ++(1cm,2cm);
\path[fill=blue, fill opacity=.1] (-2,-2) rectangle ++(1cm,2cm);
\path[fill=blue, fill opacity=.1] (.2,-2) rectangle ++(1cm,2cm);
\path[fill=blue, fill opacity=.1] (-3.2,-2.1) rectangle ++(4.5cm,2.2cm);
\node at (2.2,-1) [color=blue, opacity=.4] {{\boldmath \huge $\cdots$}};
\path[fill=blue, fill opacity=.1] (3,-2.1) rectangle ++(4.5cm,2.2cm);
\end{tikzpicture}}
}_{\mbox{$\kappa$-many blocks}}$
\caption{Subdivison of $[w]_a= \sigma_a(w)$ in the local $a$-structure
  $\M\restr[w]_a$
  for an inquisitive assignment $\Pi = \Sigma_a(w)$ that is $\kappa$-regular at
  granularity~$m$, where $\rho_m(\Pi \setminus \emptyset) = \{ \alpha_1,\alpha_2,\ldots, \alpha_n
  \}$, $\kappa = \{ 1,2, \ldots \}$ and each one of the 
  contributions $s(\alpha_j,i)$ in block $B_i$ is 
  one $\subset$-maximal element in $\Pi$, so that there are precisley
  $\kappa$ many such for each available $\simm$-type of nonempty information
  states in $\Pi$.}
\label{regularfig}
\end{figure}
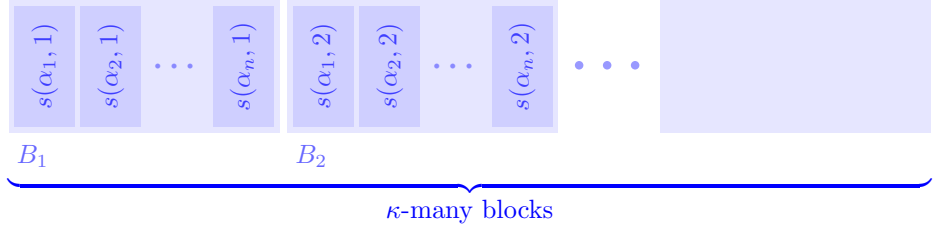

Towards normalised patterns for $\hat{\Pi} = \hat{\Sigma}_a(w) \subset \wp([w]_a)$  we think of decomposing 
$[w]_a$ into disjoint subsets $B_i \subset [w]_a$, which we call
\emph{blocks}, such that every
nonempty information state $s \in \hat{\Pi} \setminus \{
\emptyset \}$ is contained in
one of these, and such that each block $B_i$ is a
disjoint union of $\subset$-maximal information states in $\hat{\Pi}$, 
one for each
$\rho_m(s)$, $s \in \hat{\Pi} \setminus \{ \emptyset \}$.
In particular $\hat{\Pi}$ will be the union of its restrictions
$\hat{\Pi}\restr B_i := \{ s \in \hat{\Pi} \colon s \subset B_i \}$
to the $\kappa$ many blocks, as illustrated in Figure~\ref{regularfig}.

\bD
\label{regulardef}
The inquisitive assignment $\Pi = \Sigma_a(w) \subset \wp([w]_a)$
in the local structure $\M \restr [w]_a$ is \emph{$\kappa$-regular
at granularity~$m$} if 
$[w]_a = \dot{\bigcup}_{i < \kappa} B_i$ is a disjoint union of
precisely $\kappa$ many blocks $B_i \subset [w]_a$ such that
$\Pi = \bigcup_i \Pi \restr B_i$ and
\bre
\item[--]  
each of the restrictions $\Pi \restr B_i$ realises every
$\rho_m$-colour combination
$\alpha = \rho_m(s)$ for $s \in \Pi \setminus \{ \emptyset \}$ 
in precisely one $\subset$-maximal element $s(\alpha,i) \in \Pi \restr
B_i$; 
\item[--]
these maximal elements $s(\alpha,i) \in \Pi \restr
B_i$ are disjoint for distinct $\alpha$ and form a partition of $B_i$. 
\ere
\eD

It then follows that $\Pi$ itself is generated by the combination of
the downward closures of the family of disjoint maximal elements
$s(\alpha,i)$ across all $\alpha$ and $i$, which partition 
$[w]_a = \dot{\bigcup}_{\alpha,i} s(\alpha,i)$.

\bL
\label{kapparegularitylem}
Let $\M = (W,(\Sigma_a),V)$ be an inquisitive epistemic model that is
sufficiently rich in relation to given parameters $\kappa, K$ and
the cardinality of $W/{\sim}$; let $m \colon W \rightarrow \N$ be any
uniformly bounded function
prescribing finite levels of granularity for the $\sim^m$-local structures
$\M\restr [w]_a$ at $w \in W$. Then there is a globally bisimilar
companion $\hat{\M} = (W,(\hat{\Sigma}_a),V)$
such that $\hat{\M},w \sim \M,w$ for all $w \in W$
(so that the $\sim$-colouring $\rho$ on $W$ is unaffected)
that is $K$-rich and in which every local
structure $\hat{\M}\restr [w]_a$ is $\kappa$-regular at granularity $m(w)$.
\eL

\prf
It suffices to restructure every $\Pi = \Sigma_a(w) \subset \wp([w]_a)$
individually into a $\hat{\Pi}$ that is $\kappa$-regular 
w.r.t.\ $\rho_m$ for $m = m(w)$ and remains 
$K$-rich. Provided that $[w]_a$ has sufficiently many realisations
of every $\simm$-type (as guaranteed by a sufficiently high level of richness in $\M$), 
$[w]_a$ can be split into $\kappa$ many disjoint blocks each of which
is split into disjoint $\subset$-maximal new information states, 
precisely one for every $\simm$-type $\alpha = \rho(s)$ for $s \in
\Sigma_a(w)$, realising every $\simm$-type of its worlds at least $K$ times. 
\eprf

\bO
\label{richregularobs}
Every inquisitive epistemic model $\M$ admits, for $K \in \N$, finite bisimilar
coverings $\hat{\M}$ that are $K$-rich. 
As a consequence there are, for given $K,M \in \N$, globally bisimilar
inquisitive epistemic models $\hat{\M}$ sharing  
the same induced Kripke model $\str{K}(\hat{\M})$
with variations of just the inquisitive assignments, for every given
$M$-bounded granularity assignment $m \colon \hat{W} \rightarrow \N$,
that make them both $K$-rich and $\kappa$-regular 
w.r.t.\ that granularity assignment $m$.
Such $\hat{\M}$ can be chosen to be finite for finite
$\M$ and $\kappa$.
\eO

\subsubsection{Global pre-processing} 
\label{epicovsec}

The desired pre-processing needs to unclutter and smooth out the 
local overlap patterns among local structures $\str{M}\restr[w]_a$
in $\expl{\str{M}}$ in such a manner that the $\ell$-bisimulation type determines 
the first-order behaviour of $\ell$-neighbourhoods up to quantifier
rank $r$ (cf.\ Proposition~\ref{summaryGaifmanupgradeprop}).
This requires a local unfolding of the connectivity in the
underlying Kripke structures and uses ideas developed in 
\cite{OttoAPAL04} to 
eliminate incidental cycles and overlaps 
in finite bisimilar coverings --- ideas which are applied
to plain $S5$ Kripke structures in~\cite{DawarOttoAPAL09}. 
Such overlaps, which up to bisimulation equivalence are incidental (i.e.\
cannot be prescribed by the bisimulation type), occur for instance whenever an
$a$-class and an $a'$-class for $a \not= a'$ share more than one
world. These worlds would be linked by both an $R_a$-edge and an
$R_{a'}$-edge.
A bisimilar companion obtained from a classical tree unfolding
(by the appropriate closure operations that turn it into an $S5$
frame) would separate these two edges. Similarly any cycle formed by a
finite sequence of overlapping $a_i$-classes that returns to its
starting point would not be stable under appropriate bisimilar unfoldings.  

For a given epistemic model $\M$ consider the induced Kripke structure 
$\str{K} = \str{K} (\M) = (W,(R_a)_{a \in \A},V)$, which is a 
single-sorted, relational multi-modal $S5$ structure 
with equivalence relations $R_a$ derived from the
inquisitive epistemic frame $\M$, with equivalence classes
$[w]_a= \sigma_a(w)$. Generic constructions from~\cite{OttoAPAL04}, 
based on products with suitable Cayley groups, yield finite bisimilar 
coverings $\pi \colon \hat{\str{K}} \rightarrow \str{K}$ 
by $S5$ structures $\hat{\str{K}}$ that are $N$-acyclic 
in the following sense.

\begin{figure}
\centerline{%
$\begin{tikzpicture}[scale=.5]%
   \foreach \b in {1,2}
   {
    \coordinate (v\b) at ({(\b*180)+(0)}:1.1){};
    \node at (v\b) {$\bullet$};
   }
\def\elength{1.45cm}
\def\ewidth{.6cm}
\def\sshift{.15cm}
\def\lshift{1.1cm}
\path[fill=blue, fill opacity=0.2, rotate=0]
($(v1)+({\lshift},{\sshift})$) ellipse ({\elength} and {\ewidth});
\path[fill=red, fill opacity=0.2, rotate=180]
($(v2)+({\lshift},{\sshift})$) ellipse ({\elength} and {\ewidth});
\draw[red,opacity=.6] ($(v1) +(0,-0.02)$) -- ($(v2) +(0,-0.02)$);
\draw[blue,opacity=.6] ($(v1) +(0,0.02)$) -- ($(v2) +(0,0.02)$);
   \foreach \b in {1,2}
   {
    \node at (v\b) {$\bullet$};
  }
     \node at ($(v2) +(0,-3.1)$) {};
   \end{tikzpicture}$
 \qquad
\qquad
$\begin{tikzpicture}[scale=.5]%
   \foreach \b in {1,...,4}
   {
    \coordinate (v\b) at ({(\b*90)+(0)}:1.4){};
    \node at (v\b) {$\bullet$};
   }
\def\elength{1.45cm}
\def\ewidth{.6cm}
\def\sshift{.15cm}
\def\lshift{1cm}
\path[fill=blue, fill opacity=0.2, rotate=-45]
($(v1)+({\lshift},{\sshift})$) ellipse ({\elength} and {\ewidth});
\path[fill=red, fill opacity=0.2, rotate=45]
($(v2)+({\lshift},{\sshift})$) ellipse ({\elength} and {\ewidth});
\path[fill=blue, fill opacity=0.2, rotate=135]
($(v3)+({\lshift},{\sshift})$) ellipse ({\elength} and {\ewidth});
\path[fill=red, fill opacity=0.2, rotate=225]
($(v4)+({\lshift},{\sshift})$) ellipse ({\elength} and {\ewidth});
\draw[red,opacity=.6] (v1) -- (v2);
\draw[blue,opacity=.6] (v2) -- (v3);
\draw[red,opacity=.6] (v3) -- (v4);
\draw[blue,opacity=.6] (v4) -- (v1);
   \foreach \b in {1,...,4}
   {
    \node at (v\b) {$\bullet$};
  }
     \node at ($(v2) +(0,-3)$) {};
\end{tikzpicture}$   
\qquad
\qquad
$\begin{tikzpicture}[scale=.5]%
   \foreach \a in {1,...,6}
   {
       \coordinate (u\a) at ({\a*60}:2.3){};
    \node at (u\a) {$\bullet$};
   }
\def\elength{1.45cm}
\def\ewidth{.7cm}
\def\sshift{.15cm}
\def\lshift{1.2cm}
\path[fill=blue, fill opacity=0.2, rotate=60]
($(u3)+({\lshift},{\sshift})$) ellipse ({\elength} and {\ewidth});
\path[fill=red, fill opacity=0.2, rotate=120]
($(u4)+({\lshift},{\sshift})$) ellipse ({\elength} and {\ewidth});
\path[fill=blue, fill opacity=0.2, rotate=180]
($(u5)+({\lshift},{\sshift})$) ellipse ({\elength} and {\ewidth});
\path[fill=red, fill opacity=0.2, rotate=240]
($(u6)+({\lshift},{\sshift})$) ellipse ({\elength} and {\ewidth});
\path[fill=blue, fill opacity=0.2, rotate=300]
($(u1)+({\lshift},{\sshift})$) ellipse ({\elength} and {\ewidth});
\path[fill=red, fill opacity=0.2, rotate=0]
($(u2)+({\lshift},{\sshift})$) ellipse ({\elength} and {\ewidth});
\draw[red,opacity=.6] (u1) -- (u2);
\draw[blue,opacity=.6] (u6) -- (u1);
\draw[red,opacity=.6] (u5) -- (u6);
\draw[blue,opacity=.6] (u4) -- (u5);
\draw[red,opacity=.6] (u3) -- (u4);
\draw[blue,opacity=.6] (u2) -- (u3);
   \foreach \a in {1,...,6}
   {
    \node at (u\a) {$\bullet$};
   }
 \end{tikzpicture}$}

\vspace{.3cm}
\centerline{%
$\begin{tikzpicture}[scale=.5]%
   \foreach \b in {1,...,3}
   {
    \coordinate (v\b) at ({(\b*120)+(30)}:1.3){};
    \node at (v\b) {$\bullet$};
   }
\def\elength{1.45cm}
\def\ewidth{.6cm}
\def\sshift{.15cm}
\def\lshift{1.1cm}
\path[fill=blue, fill opacity=0.2, rotate=0]
($(v1)+({\lshift},{\sshift})$) ellipse ({\elength} and {\ewidth});
\path[fill=black, fill opacity=0.2, rotate=120]
($(v2)+({\lshift},{\sshift})$) ellipse ({\elength} and {\ewidth});
\path[fill=red, fill opacity=0.2, rotate=240]
($(v3)+({\lshift},{\sshift})$) ellipse ({\elength} and {\ewidth});
\draw[red,opacity=.6] (v2) -- (v3);
\draw[black,opacity=.6] (v1) -- (v2);
\draw[blue,opacity=.6] (v3) -- (v1);
   \foreach \b in {1,...,3}
   {
    \node at (v\b) {$\bullet$};
  }
     \node at ($(v3) +(0,-3.5)$) {};
\end{tikzpicture}$
\qquad\qquad
$\begin{tikzpicture}[scale=.5]%
   \foreach \a in {1,...,6}
   {
       \coordinate (u\a) at ({\a*60}:2.3){};
    \node at (u\a) {$\bullet$};
   }
\def\elength{1.45cm}
\def\ewidth{.7cm}
\def\sshift{.15cm}
\def\lshift{1.2cm}
\path[fill=blue, fill opacity=0.2, rotate=60]
($(u3)+({\lshift},{\sshift})$) ellipse ({\elength} and {\ewidth});
\path[fill=black, fill opacity=0.2, rotate=120]
($(u4)+({\lshift},{\sshift})$) ellipse ({\elength} and {\ewidth});
\path[fill=red, fill opacity=0.2, rotate=180]
($(u5)+({\lshift},{\sshift})$) ellipse ({\elength} and {\ewidth});
\path[fill=blue, fill opacity=0.2, rotate=240]
($(u6)+({\lshift},{\sshift})$) ellipse ({\elength} and {\ewidth});
\path[fill=black, fill opacity=0.2, rotate=300]
($(u1)+({\lshift},{\sshift})$) ellipse ({\elength} and {\ewidth});
\path[fill=red, fill opacity=0.2, rotate=0]
($(u2)+({\lshift},{\sshift})$) ellipse ({\elength} and {\ewidth});
\draw[red,opacity=.6] (u1) -- (u2);
\draw[black,opacity=.6] (u6) -- (u1);
\draw[blue,opacity=.6] (u5) -- (u6);
\draw[red,opacity=.6] (u4) -- (u5);
\draw[black,opacity=.6] (u3) -- (u4);
\draw[blue,opacity=.6] (u2) -- (u3);
   \foreach \a in {1,...,6}
   {
    \node at (u\a) {$\bullet$};
   }
 \end{tikzpicture}$ \qquad\qquad\nt}
\caption{Example of particularly simple unfoldings in bisimilar coverings of 
  $S5$ Kripke frames for agents \emph{red}, \emph{blue},
  \emph{black}; top row: unfolding of a $2$-cycle (overlap of two
  equivalence classes in two worlds) into a $4$- or
  $6$-cycle;
  below: unfolding of a $3$-cycle into a $6$-cycle
  (just non-trivial equivalence classes colour-coded and worlds in
  overlaps and
  links in induced cycles marked).}
\label{Nacycfig}
\end{figure}
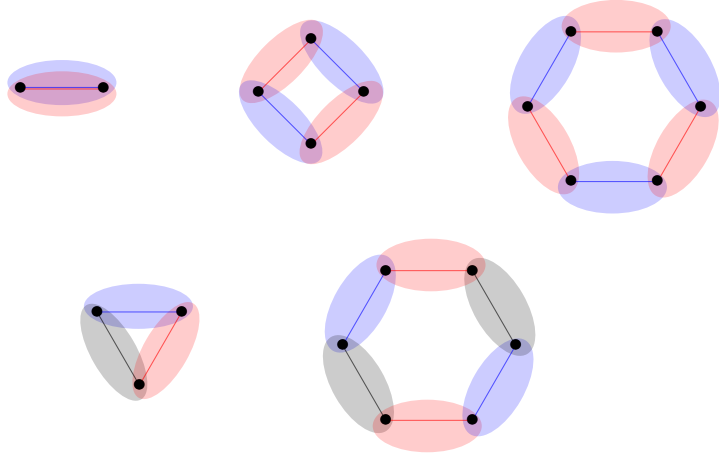

\bD
\label{Nacycdef}
  For $N \geq 2$,
  an epistemic model $\M$ is called \emph{$N$-acyclic} if
  the induced $S5$ Kripke frame $(W,(R_a)_{a \in \A})$ underlying 
  $\str{K} = \str{K} (\M)$ is $N$-acyclic in the sense that 
  there is no non-trivial cyclic pattern of 
  overlapping $a_i$-classes of length up to $N$: for $2 \leq n \leq N$, there
  is no cyclically indexed family $([w_i]_{a_i})_{i \in \Z_n}$ such that
   $[w_i]_{a_i} = [w_{i+1}]_{a_i}$ while $w_{i+1} \not= w_i$ and
   $a_{i+1} \not= a_i$. 
\eD

This acyclicity condition implies in particular that $a_i$-classes
and $a_j$-classes for agents $a_i \not= a_j$ overlap in at most a
single world: if $[w]_{a_1} \cap [w]_{a_2} \not= \{ w \}$ the two
classes would form a non-trivial $2$-cycle.
For $N \geq 2 \ell$ it moreover implies that
the overlap pattern between $a$-local structures is simply
tree-like within Gaifman radius~$\ell$: any two
distinct short sequences of such overlaps
emanating from the same world cannot meet elsewhere.

By~\cite{OttoAPAL04}, this form of $N$-acyclicity can be achieved, for
any given $N$, in finite bisimilar coverings in which every $a$-class is 
bijectively related to an $a$-class in the original structure
by the projection map of the covering, which induces the bisimulation
(cf.~Figure~6 for the basic intuition). 
We briefly outline the extension of such bisimilar 
coverings of $\str{K}(\M)$ to bisimilar coverings
of the inquisitive epistemic model $\M$. 
Let $\pi \colon \hat{\str{K}} \rightarrow \str{K}$
be such an $N$-acyclic finite bisimilar covering, which in
particular means that the graph of the covering projection $\pi$ is a
bisimulation relation. Since $a$-classes in $\hat{\str{K}}$
are $\pi$-related isomorphic copies of corresponding $a$-classes in
$\str{K}$, we may consistently endow them with an
inquisitive $\Sigma_a$-assignment as pulled back from $\M$,
to obtain a natural derived finite inquisitive bisimilar covering
\[
\pi \colon \hat{\M}
\longrightarrow \M = (W,(\Sigma_a),V)
\]
in the sense of Definition~\ref{inqbisimcoverdef}.
The inquisitive epistemic model $\hat{\M} = (\hat{W},(\hat{\Sigma}_a),
\hat{V})$ is built on the basis of 
the covering $S5$ Kripke structure $\hat{\str{K}}$ 
so that $\str{K}(\hat{\M}) =\hat{\str{K}}$. More specifically, 
for $\hat{w} \in \hat{W}$ consider its $\hat{R}_a$-equivalence class
$[\hat{w}]_a\subset \hat{W}$, which is bijectively related by
$\pi$ to the $R_a$-equivalence class of $w = \pi(\hat{w})$ in
$\str{K}$, $[w]_a \subset W$. In restriction to  
$[\hat{w}]_a$ we put 
\[
\barr{rcl}
\hat{\Sigma}_a \restr [\hat{w}]_a\colon [\hat{w}]_a
&\longrightarrow & 
\wp\wp( [\hat{w}]_a)
\\
\hnt
\hat{u} &\longmapsto&
\{ 
\pi^{-1}(s) \colon s \in \Sigma_a(w) \subset \wp( [w]_a)
\},
\earr
\]
where $\pi^{-1}$ refers to the inverse of the local bijection 
between $[\hat{w}]_a $ and $[w]_a$ induced by $\pi$.
In other words, we make the following diagram commute and note that 
the relevant restriction of $\pi$ in the right-hand part of the
diagram is bijective and part of a
bisimulation at the level of the underlying $S5$ Kripke structures, which
means in particular that it is compatible with the given propositional
assignments:
\[
\xymatrix{
*++{\hat{\M},\hat{w}} \ar[dd]_{\pi } & 
[\hat{w}]_a
\ar[rr]^{{\hat{\Sigma}_a}} \ar[dd]_{\pi/\pi^{-1}} 
&&
\wp\wp([\hat{w}]_a)
\ar[dd]^{\pi/\pi^{-1}} 
\\
\\
*++{\M,w} & 
[w]_a 
\ar[rr]_{{\Sigma_a}}
\ar[uu]
&& 
\wp\wp([w]_a)
\ar[uu]
}
\]

We merely need to check that the resulting model $\hat{\M}$ is again an inquisitive
epistemic model, and that $\pi$ induces an inquisitive bisimulation, not
just a bisimulation at the level of the underlying single-sorted,
multi-modal $S5$-structures $\hat{\str{K}} = \str{K}(\hat{\M})$ and 
$\str{K} = \str{K}(\M)$. 
That $\hat{\M}$ is an inquisitive epistemic model is straightforward 
from its construction: 
$\hat{\Sigma}_a$-values are constant on the $a$-classes induced
by $\hat{R}_a$ just as 
$\Sigma_a$-values are constant across $a$-classes induced by
$R_a$; and $\hat{w} \in \hat{\sigma}_a(\hat{w})$ follows from the fact that 
$w \in \sigma_a(w) = [w]_a$.

A non-deterministic winning strategy in the unbounded inquisitive
bisimulation on $\hat{\M},\hat{w}$ and $\M,\pi(\hat{w})$ is given by
the following bisimulation relation, to be viewed as a set of winning
positions for player~\PII\ in the bisimulation game:
\[
\barr{r@{\;}l}
Z = & \{ (\hat{u}, \pi(\hat{u})) \in \hat{W} \times W \colon 
\hat{u} \in \hat{W} \}
\\
\hnt
& \cup \;
\{ (\hat{s}, \pi(\hat{s}) ) \in \wp(\hat{W}) \times \wp(W) 
\colon a \in \A, \hat{w} \in \hat{W},
\hat{s} \in \hat{\Sigma}_a(\hat{w}) \},
\earr
\]
which is the natural lifting 
of $\pi$ to information states. It is easily verified that
player~\PII\ can stay in $Z$ in response to any challenges by
player~\PI\ ($Z$ satisfies the usual back\&forth conditions).
We summarise these findings as follows.

\bL
\label{S5unfoldlem}
Any inquisitive epistemic model $\M = 
(W, (\Sigma_a),V)$ admits, for every $N \in
\N$, a finite bisimilar covering of the form
\[
\pi \colon \hat{\M} \longrightarrow \M
\]
by an inquisitive epistemic model
$\hat{\M}= (\hat{W}, (\hat{\Sigma}_a),\hat{V})$
such that
\bre
\item
$\hat{\M}$ is $N$-acyclic in the sense of Definition~\ref{Nacycdef}; 
\item
the global bisimulation induced by $\pi$ is an isomorphism 
in restriction to each local structure $\hat{\M}\restr
[\hat{w}]_a$ of $\hat{\M}$, which is isomorphically mapped by $\pi$
onto the local $a$-structure $\M\restr [\pi(\hat{w})]_a$ of $\M$. 
\ere
\eL

It is noteworthy that condition~\emph{(ii)} guarantees that 
degrees of regularity and richness are preserved in the covering, simply
because they are properties of the local structures
$\M\restr[w]_a$, which up to isomorphism are the same 
in $\M$ and in $\hat{\M}$. So, in combination 
with Observation~\ref{richregularobs},  we obtain the following.

\bP
\label{bisimpartnerprop}
For any desired values of $K,N \in \N$ and (finite or infinite)
cardinality $\kappa$, and for any finite bound $M$ on relevant
granularities, 
any world- or state-pointed inquisitive epistemic model $\M$ admits
a globally bisimilar $N$-acyclic companions $\hat{\M} \sim \M$,
all based on the same $N$-acyclic induced Kripke model $\str{K}(\hat{\M})$,
whose epistemic assignments can be chosen $K$-rich and $\kappa$-regular
w.r.t.\ any given uniformly bounded
granularity assignment $m \colon \hat{W}\rightarrow \N$.
Such $\hat{\M}$ can be chosen finite for finite $\M$
and $\kappa$. In the non-trivially state-pointed case of $\M,s$,
companions $\M_s$ based on the dummy agent~$s$ admit
bisimilar companions $\hat{\M},\hat{s} \sim \M_s, s$
with these same properties. 
\eP

The non-trivially state-pointed case may deserve further
comment. The desired $\hat{\M}$ in this case is 
based on the augmentation $\M_s$ with the extra agent~$s$.
The local pre-processing produces in particular a rich and
$\kappa$-regular
$s$-local structure over a set of worlds~$\hat{s}$; $N$-acyclicity
for $N > 2\ell$ in the global pre-processing  then serves to
scatter the worlds in $\hat{s}$ within $\str{K}(\hat{\M})$ w.r.t.\ $\A$ in the sense that no
two $(E_a)_{a \in \A}$-paths of length up to~$\ell$ from distinct worlds in
$\hat{s}$ can meet: they are directly $E_s$-linked and $N$-cyclic
patterns are ruled out.

\subsection{Proof of the characterisation theorem} 
\label{epicharsec}
In the remainder of this section we complete the proof of our main theorem, Theorem~\ref{main2}, 
the precise statement of which was given
in Section~\ref{bisimconstraintsec}
on page~\pageref{th:restated}. 

The proof of this theorem has, through
preparations in the preceding sections, already been reduced to the upgrading target in
equation~(\ref{Gaifmanequiveqn}). For the world-pointed case:
\[
  \expl{\hat{\str{M}},\hat{w}}
\;\equiv_{r}^{\ssc (\ell)}\;
\expl{\hat{\str{M}}',\hat{w}'},
\]
where we now take $\hat{\str{M}}$ and $\hat{\str{M}}'$ to be locally full
relational encodings of bisimilar
companions of any given pair of $\simn$-bisimilar pointed models
$\M(\str{M})$ and $\M(\str{M}')$, which are obtained through the detour
summarised in Proposition~\ref{bisimpartnerprop}, for suitable values
of the parameters of $K$ and~$N$ for richness and acyclicity, and
$\kappa$ and $m$ for regularity. 
The parameters of the target equivalence, $\ell$ and~$r$, are
determined by the quantifier rank of the given $\sim$-invariant $\phi
\in\FO$ that is to be shown to be $\simn$-invariant, for a value
of~$n$ still to be determined by the Ehrenfeucht--Fra\"\i ss\'e 
argument for~$\equiv^{\ssc (\ell)}_r$.

\medskip
Recall from our preliminary discussion at the end
Section~\ref{FOEFGaifsec}
how the desired upgrading argument for the crucial compactness claim in
line with Proposition~\ref{summaryGaifmanupgradeprop},
calls for a strategy in the $r$-round first-order
Ehrenfeucht--Fra\"\i ss\'e\ game over the $\ell$-neighbourhoods
of exploded views of the pre-processed models. Due to the nature of
the exploded views as conglomerates of interconnected $a$-local
structures, the desired strategy analysis focuses on a separation between
\bae
\item
the global analysis of the $\FO$-game in the $\ell$-neighbourhoods, now 
involving a tree-like connectivity pattern between $a$-local component
structures in sufficiently acyclic models.
\item
the $a$-local concerns, now involving sufficiently rich and suitably
regular models, where the $\FO$-game over these two-sorted structures can be
analysed in terms an $\MSO$-analysis over their first sorts, to be
treated first, in Section~\ref{MSOequivsec}. 
\eae

The granularity for the regularity requirements in~(b) will depend on the
distance from the distinguished world or state, since control over
finite levels of bisimulation equivalence decreases with this
distance from the root in the tree-like core structure. 
The resulting strategy advice for player~\PII\ is brought together
in Section~\ref{FOequivsec}.
We first turn to these $a$-local component structures.

\subsubsection{MSO-equivalence over local structures}
\label{MSOequivsec}
Recall Definition~\ref{localastrgranulardef} of local structures
$\M\restr[w]_a$ with colouring $\rho_m$ that marks out 
$\simm$-types in the surrounding $\M$ or $\M'$.
In the locally full relational encoding $\str{M} = \str{M}^\lf(\M)$,
we correspondingly consider
$\str{M}\restr[w]_a,\rho_m$ as a mono-modal relational model.
We apply Ehrenfeucht--Fra\"\i ss\'e techniques for $\MSO$ in order to
obtain guarantees for levels of first-order equivalence $\equiv_k$ over these
two-sorted local $\simm$-structures, assuming a
suitable degree of richness and regularity w.r.t.\ $\rho_m$. 
Richness matters since $\FO$
can assess cardinality distinctions between information states
in the second sort, but just up to certain finite thresholds.
These thresholds are inherited from $\MSO$-definability over
plain coloured sets over the first sort. 

To compare the sizes of sets $P,P'$ up to some
threshold~$d \in \N$, we introduce this cardinality equality with
cut-off: 
\[
  |P| =_d |P'| \quad \colon\!\!\Leftrightarrow \quad  |P|= |P'| \mbox{ or } |P|,|P'| \geq d.
 \]

For the analysis of 
a set $W$ with any $C$-colouring $\rho \colon W \rightarrow C$,
we also use its representation as a
plain relational structure $(W,(W_c)_{c \in C})$ with colour
classes~$W_c$ viewed as the interpretations of unary relations
(monadic predicates).
So we interchangeably think of representations of $W,\rho$ as a set with a
colouring $\rho \colon W \rightarrow C$ or as a relational structure
partitioned into these colour classes. In our applications, $W$ will
always be the first sort $[w]_a$ of some local $\simm$-structure, 
coloured by $\rho_m$ with the finite colour set for $\simm$-types of
worlds.

With any subset $P \subset W$ we correspondingly associate the family
\[
  (P_c)_{c \in C} \;\mbox{ where }\; P_c := \{ w \in P \colon \rho(w) = c
  \} =
  P \cap W_c.
\]

For the comparison of $P \subset W$ with $P' \subset W'$ from
two $C$-coloured sets $W,\rho$ and $W',\rho'$ we define the equivalence
relation $\approx^C_d$ according to 
\[
P \approx^C_d P'  \quad\mbox{ iff }\quad |P_c| =_d |P_c'| \mbox{ for all } c \in C;
\]
and generalise this to $k$-tuples of subsets 
$\Pbar \in (\wp(W))^k$ 
and
$\Pbar' \in (\wp(W'))^k$  
by putting
\[
\Pbar \approx^C_d \Pbar' \quad\mbox{ iff }\quad 
\zeta(\Pbar) \approx^C_d \zeta(\Pbar')
\mbox{ for all boolean terms }\zeta(\Xbar). \footnotemark
\]
\footnotetext{Here a \emph{boolean term} is any term in the functional language of 
Boolean algebras, with binary operations for union and intersection and a
unary operation for complementation.}%

The following is folklore, but 
we indicate the straightforward proof below.

\bL
\label{EFMSOlem}
For $C$-coloured sets $W,W'$
with $k$-tuples $\Pbar \in \wp(W)^k$ and
$\Pbar' \in \wp(W')^k$
of subsets, 
and for $d= 2^q$:
\[
\Pbar \approx^C_d \Pbar'
\quad
\Rightarrow
\quad 
(W,(W_c)_{c \in C},\Pbar) \equiv^{\MSO}_q (W',(W_c')_{c \in C},\Pbar').
\]
\eL

Note that assignments to first-order element variables $x$ are
implicit in the form of singleton sets for $\{ x\}$ (which induces
but a small constant shift in quantifier rank levels).

\prf[Proofsketch.]
The proof is by induction on $q$, and we first look at the case with
trivial colouring (all elements are of the same colour).
For the induction step assume that $\Pbar \approx_{2d}^C \Pbar'$
and suppose w.l.o.g.\ that the first player proposes 
a subset $P \subset W$ so that the second player needs to find a response 
$P' \subset W'$ such that $\Pbar P \approx^C_d \Pbar' P'$.
Decompose $P$ and its complement $\bar{P}$ 
into their intersections with the atoms of the 
boolean algebra generated by $\Pbar$ in $\wp(W)$. 
Then each part of these partitions of $P$ and $\bar{P}$
can be matched with a subset of the corresponding atoms of the
boolean algebra generated by $\Pbar'$ in $\wp(W')$ in such a
manner that corresponding parts match in the sense of $=_d$.
This just uses the assumption that $|\zeta(\Pbar)| =_{2d}
|\zeta(\Pbar')|$ for each $\zeta$. In the case of non-trivial
$C$-colourings, the same argument can be played out within each pair
of colour classes $W_c$ and $W_c'$ after decomposition of $P$ and
$\bar{P}$ into their $C$-partitions, and recombination of $P'$ as the
union of the matching 
$P_c'$ obtained for the $P_c$.
\eprf

We next transfer this $\MSO$-equivalence between
$\rho_m$-coloured first sorts of suitable pre-processed 
local component structures to $\FO$-equivalence between
the two-sorted relational encodings of these 
local $\simm$-structures $\str{M}\restr[w]_a,\rho_m$ and
$\str{M}'\restr[w]_a,\rho_m$ that locally realise the same
$\simm$-types of information states and
are $\kappa$-regular at this granularity~$m$ for the same value of $\kappa$.

Towards an interpretation of the $2$-sorted full
relational encodings of local $\simm$-structures like
$\str{M}\restr[w]_a,\rho_m$ in terms of $\MSO$-definability over $[w]_a$,
it is necessary to make available the
information states in $\Pi:=\Sigma_a(w)$. Due to
$\kappa$-regularity at granularity~$m$, this reduces to marking out
the $\subset$-maximal informatiobn states $s(\alpha,i)$ in each block
$B_i$, for every one of the finitely many colour combinations $\alpha =
\rho_m(s) \in \rho_m(\Pi) \setminus \{ \emptyset \}$ for non-trivial
information states~$s$ realised in
$\str{M}\restr[w]_a,\rho_m$. To this end we look at
$\MSO$-interpretability over derived single-sorted structures over
$[w]_a$. These encode the block decomposition (by an equivalence
relation~$B$)
and the $\subset$-maximal
constituents across the blocks (by finitely many unary predicates $P_\alpha$):
\[
\str{b}_m([w]_a)
\;:=\; \barr[t]{@{}l@{}} ([w]_a,\rho_m, B, (P_\alpha) )
\\
\hnt
    \mbox{ for }   \barr[t]{@{}l@{\;:=\;}l@{}} \; B & \bigcup_{i < \kappa} B_i^2 = \{ (u,v) \in [w]_a^2 \colon u,v
\mbox{ in the same block}\},
\\
\hnt
P_\alpha & \bigcup_{i < \kappa} s(\alpha,i).
\earr
\earr
\]

The full relational encoding of a
$\simm$-local structure $\str{M}\restr [w]_a,\rho_m$ becomes
$\MSO$-interpretable over these $\str{b}_m([w]_a)$. In particular,  for
  $s \in \wp([w]_a)$ we have $s \in \Sigma_a(w)$ if, and only if 
  \[
   \textstyle
    \str{b}_m([w]_a),s \models   \exists X \Bigl(
    s \subset X \wedge 
    \forall x \forall y \bigl( (Xx \wedge Xy) \rightarrow Bxy \bigr)
    \wedge
    \bigvee_\alpha X \subset P_\alpha\Bigr).
    \]

It follows that, for some fixed set-off $c$
stemming from the underlying interpretation, and tuples
$\sbar,\sbar'$ of set parameters over $[w]_a$ and $[w']_a$,
\begin{equation}
\str{b}_m([w]_a), \sbar
\equiv^\MSO_{r+c}\str{b}_m([w']_a), \sbar'
\;\; \Rightarrow \;\;
\str{M}\restr [w]_a,\rho_m, \sbar
\equiv^\FO_{r}
\str{M}'\restr [w']_a,\rho_m,\sbar'
\label{MSOtoFOeqn}
\end{equation}
    
In the analysis of $\FO$-equivalences between
local $\simm$-structures (as part of the Ehrenfeucht--Fra\"\i ss\'e
argument over exploded views as in the upgrading as in
Figure~\ref{S5extendedupgradefigure})
we shall encounter choices of elements inside these $2$-sorted
structures as parameters. These may be worlds $v \in [w]_a$ (first
sort) or subsets $s \subset [w]_a$ (second sort). Looking at these
in the $\MSO$-perspective of the above interpretation, we can treat
both kinds as monadic second-order parameters by associating
$v \in [w]_a$ with $\{ v \} \subset [w]_a$. In particular, a source 
world $w \in [w]_a$ for the local $\simm$-structure associated with
$[w]_a$ can be fed as a distinguished world in $\str{M}\restr [w]_a,\rho_m$
and as a singleton set parameter $\{ w\} \subset [w]_a$
in the $\MSO$-analysis of $\str{b}_m([w]_a)$. 
    
\bL
\label{MSOequiv}
  Any pair of $m$-bisimilar pointed $\simm$-local structures
  $\M\restr [w]_a,w$ and $\M'\restr [w']_a,w'$ for $m \geq 1$
  that are both $\kappa$-regular for the same
  (finite or infinite) value of~$\kappa$, realise the
  same colour combinations $\alpha \subset \rho_m([w]_a) = \rho_m([w']_a$
  as nonempty contributions to the partitions of
  $[w]_a = \dot{\bigcup}_\alpha P_\alpha$ and
  $[w']_a = \dot{\bigcup}_\alpha P_\alpha'$, respectively.
  Moreover 
 \[
    \str{b}_m([w]_a), \{w\}
    \equiv^\MSO_q
   \str{b}_m([w']_a), \{w'\},
 \]
 provided both are sufficiently rich in relation to~$q$.
 \eL

 It then follows from equation~\ref{MSOtoFOeqn} that in this situation
 also
 \[
   \str{M}\restr [w]_a,\rho_m,w  \equiv^\FO_k
   \str{M}'\restr [w']_a,\rho_m,w',
 \]
given the appropriate
level of richness for~$q = k+c$ (cf.~Corollary~\ref{localupgradecor} below).

\prf
Note that in the degenerate case of $\kappa =1$ (regularity in a single block,
trivial $B$), the claim follows directly from Lemma~\ref{EFMSOlem}. 
The extension to finite $\kappa$ then follows from compatibility of
$\equiv^\MSO$ with disjoint unions of structures (in application to
matched blocks $B_i$ and $B_i'$ as marked out by the equivalence
relations $B$ and $B'$). But also in the general case of
possibly infinite $\kappa$, any choice of a bijection between
the $\kappa$ many blocks of either side, starting with the blocks $B_0$ and
$B_0'$ that contain the distinguished worlds $w$ and $w'$,
respectively, allows us to put together a response to a challenge by the
first player from responses in the restrictions of the games over
matched pairs of blocks, just as in the proof sketch for
Lemma~\ref{EFMSOlem}, but w.r.t.\ a decomposition into a possibly
infinite number of non-trivial restrictions to individual blocks.
\eprf

In the comparisons between $\kappa$-regular local $\simm$-structures 
$\M \restr [w]_a, \rho_m$ and $\M' \restr [w']_a, \rho_m$ it is
essential that $m \geq 1$. As remarked in the lemma, 
$\sim^1$  between any two distinguished worlds,
$\M \restr [w]_a, \rho_m, w \sim^1 \M' \restr [w']_a, \rho_m,w'$,
implies that these two local structures are globally bisimilar and in particular
realise the same $\simm$-types both of worlds $v \in [w]_a$ and $v'
\in [w']_a$ and of information states
$s \in \Pi := \Sigma_a(w)$ and $s \in \Pi' := \Sigma_a'(w)$.
The encodings of the block structure and $\subset$-maximal information
states in the $\kappa$-regular inquisitive assignments $\Pi$ and $\Pi'$ were
essential for the background $\MSO$-analysis over the first sorts, but 
we use these results towards the $\FO$-analysis over their $2$-sorted
full relational encodings. In these terms we have established the
following.

\bC
\label{localupgradecor}
If $\M$ and $\M'$ are both $\kappa$-regular (for the same $\kappa$)
at granularity~$m$ and sufficiently rich (in
relation to $k \in \N$), then the following holds for any two worlds
$v \in [w]_a$ and $v' \in [w']_a$:
\[
  \barr{r@{\;\;\Rightarrow\;\;}l}
  \M\restr[w]_a, \rho_m, v \sim^1 \M'\restr[w']_a,\rho_m,v'
  &
  \\
  \hnt
  \M\restr[w]_a, \rho_m, v \sim \M'\restr[w']_a,\rho_m,v'
  &
 \str{M}\restr [w]_a, \rho_m, v \equiv_k \str{M}\restr [w']_a,
 \rho_m ,v'.
 \earr
\]
\eC

So Observation~\ref{richregularobs} and
Proposition~\ref{bisimpartnerprop}, 
any two inquisitive epistemic models
$\M$ and $\M'$ admit, for any target value $N, k \in \N$,
globally bisimilar $N$-acyclic companions $\hat{\M} \sim \M$ and 
$\hat{\M}' \sim \M'$, finite if $\M$ and $\M'$ are finite,
whose locally full relational encodings satisfy
the above transfer within every one of their local
$\simm$-structures. 
This conditional guarantee of $\FO$-indistinguishability
is good for any allocation of the granularities $m \geq 1$, which
moreover may take different values for different pairs of local structures that are to
be matched. It is important that the local analysis reflected in the
corollary obviously imports -- through $\rho_m$ --
global information about $\simm$-types in
$\hat{\M}$ and $\hat{\M}'$ into the local structures at hand, but that
it does not rely on a uniform allocation of granularity across $\hat{\M}$ or $\hat{\M}'$. 

\subsubsection{Local FO-equivalence of exploded views}
\label{FOequivsec}
Returning to the global view of the underlying models
$\M =  (W,(\Sigma_a)_{a \in \A},V)$ and $\M' =
(W',(\Sigma_a')_{a \in \A},V')$ and their exploded views, and assuming
$\M,w \;\simn\; \M',w'$, 
we recall equation~(\ref{Gaifmanequiveqn}) and its r\^ole on the
lefthand side of Figure~\ref{S5extendedupgradefigure} for the
world-pointed case. We first cover this basic, world-pointed version.
The adaptation to the non-trivially state-pointed variant, involving the dummy agent
$s$, will then be straightforward. 
For the world-pointed version, 
we need to lift the $[\,\cdot\,]_a$-local equivalences
(at suitable granularities
$1 \leq m\leq n$) in matching local
$\simm$-structures
to the $\ell$-local equivalence 
\[
  \expl{\hat{\str{M}},\hat{w}}
\;\equiv_{r}^{\ssc (\ell)}\;
\expl{\hat{\str{M}}',\hat{w}'}
\]
around the distinguished worlds $\hat{w}$ and $\hat{w}'$.
We aim to employ the criterion of $N$-acyclicity from
Definition~\ref{Nacycdef}, more specifically its consequences for
connectivity patterns in the exploded views. In order to reduce the
notational overhead in the following we assume w.l.o.g.\ that
$\hat{\str{M}} = \str{M}$ and $\hat{\str{M}}' = \str{M}'$, or that  
$\M$ and $\M'$ themselves have been pre-processed in the sense of
Proposition~\ref{bisimpartnerprop}. In the exploded view of such
$N$-acyclic $\M$, any two local structures
 $\str{M}\restr[w_1]_{a_1} \not= \str{M}\restr[w_2]_{a_2}$ can be linked through at
  most one element $v$ shared by $[w_1]_{a_1}$ and $[w_2]_{a_2}$. This
  follows from Definition~\ref{Nacycdef} for any $N \geq2$.
  Such a link, if it exists, relates the copies of $v$ in the two component
  structures of $\expl{\str{M}}$ by  a path of length $2$, which is formed by
  two $R_I$-edges emanating from the node $v \in W$ in the core part.
  So this core part $(W,(E_a)_{a \in \A})$, which is the underlying $S5$
  frame, provides the only links between individual local constituents
  $\str{M}\restr[w_i]_{a_i}$ in the exploded view $\expl{\M}$. 
  For $N$-acyclic $\M$, this core part is $N$-acyclic, so 
  that there cannot be any non-trivial overlaps 
  between any two distinct threads  of lengths up to $\ell$
  of local constituents $\str{M}\restr[w_i]_{a_i}$
  in the exploded view $\expl{\M}$ emanating from the same
  $v \in W$, provided $N > 2\ell$.

 \bL
\label{inqS5upgradelem}
Let $\M$ and $\M'$ be $N$-acyclic for 
$N > 2\ell$, $\kappa$-regular (for the same $\kappa$)
at granularity~$m$ (depending on depth in the tree-like exploded
views, i.e.\ on distance from the distinguished worlds~$w,w'$)
and sufficiently rich to support the claim of 
Corollary~\ref{localupgradecor}. 
Then their exploded views
$\expl{\M} = \expl{\str{M}^\lf(\M)}$ and 
$\expl{\M'} = \expl{\str{M}^\lf(\M')}$
satisfy the following for any pair of worlds such that
$\M,w \sim^n \M',w'$ for $n \geq \ell+1$:
\[
  \expl{\M},w \equiv_r^{\ssc (\ell)} \expl{\M'},w'
  \;\mbox{ i.e.\ }\;
 \expl{\M} \restr N^\ell(w),w \equiv_r \expl{\M'} \restr N^\ell(w'),w'.
\]
\eL

The proof of this technical core lemma is given in terms of 
a winning strategy for player~\PII\ in the $r$-round game 
on the pointed $\ell$-neighbourhoods of 
$w$ and $w'$ in $\expl{\M}$ and $\expl{\M'}$.
As outlined at the start of this section, Section~\ref{epicharsec},
such a strategy
can be based on a combination of strategies in 
\bae
\item
the $r$-round first-order game on the tree-like underlying 
$S5$ Kripke structures on the first sort (worlds) with $a$-classes for accessibility
and a colouring of worlds that reflects sufficient finite
approximations to their bisimulation types in $\M$,  
which in particular suffice to guarantee $1$-bi\-si\-mu\-la\-tion 
equivalence and hence also $\equiv_r$ of local component structures 
of appropriate granularities;%
\footnote{It is important here that just certain finite levels of bisimilarity, not the full
$\sim$-type, can be controlled.}
\item
local strategies in the $r$-round games on the relational encodings of
the local component structures (with colourings 
at the appropriate level of granularity) 
induced by pebbled pairs of worlds or states.
These strategies are ultimately based on  
local $\MSO$-games over their first sorts, as discussed in connection
with Corollary~\ref{localupgradecor}.
\eae

The requirements w.r.t.\ $\sim^m$-matches
will decrease, in terms of the parameter $m$,
with distance from the source worlds $w$ and $w'$, but
keeping $m \geq 1$ to uniformly support~(a).

\prf[Proof of Lemma~\ref{inqS5upgradelem}.]
Let $\str{N} := \expl{\M} \restr N^\ell(w)$ and 
$\str{N}' := \expl{\M'} \restr N^\ell(w')$; for the local component
structures we directly refer to the locally full relational encodings
$\str{M} := \str{M}^\lf(\M)$ and  $\str{M}' := \str{M}^\lf(\M')$.
Regularity in matching local structures $\str{M}\restr[u]_a$ and
$\str{M}'\restr[u']_a$ is imposed at granularity $m(u) := \ell - d(w,u)
+1 = \ell - d(w',u') = m(u')$ where $d$ is the depth of $u,u'$ in the tree-like
core parts of the exploded views, i.e.\ distance from $w,w'$.

It follows that  $1 \leq m(u),m(u') \leq \ell + 1 \leq n$ is uniformly bounded
by the value~$n$ in the initial assumption about $\M,w \sim^n \M',w'$.

The strategy analysis in support of the lemma is provided in three steps:
\bre
\item[(I)] parameterisation of game positions player~\PII\ may
  have to deal with;
\item[(II)] abstraction of an invariant capturing the relevant
  information in those;
  \item[(III)] showing how player~\PII\ can maintain this invariant
    through $r$ rounds.
\ere

\paragraph*{(I) Parameterisation of game positions.} 
We consider any position $\ubar;\ubar'$ in the $r$-round first-order game 
on $\str{N}$ and $\str{N}'$, 
starting with pebbles on $w$ and $w'$ which are equivalent w.r.t.\
$\sim^{\ell+1}$ in $\M$ and $\M'$. It will be of the
form $\ubar = (u_0,\ldots,u_k)$ and 
$\ubar' = (u'_0,\ldots,u'_k)$ for some $k \leq r$ with $u_0 = w$, 
$u'_0 = w'$, and,  for $1 \leq i \leq k$, either
\bre
\item
$u_i= w_i$ and $u_i'= w_i'$ are elements of the first sort
in the core parts that are linked to $w$ and $w'$, respectively, by unique
paths of overlapping $a$-classes, or 
\item
$u_i= (w_i)_a$ and $u_i'= (w_i')_a$ are elements of the first sort
in local constituents $[w_i]_a$ and $[w_i']_a$,
for some $a = a_i \in \A$ and$w_i,w_i'$ in the core parts as in (i),
\\
or
\item 
$u_i = s_i$ and $u_i'= s_i'$ are of the second
sort (information states) in some $\wp([w_i]_a)$, respectively
$\wp([w_i']_a)$,
for some $a = a_i \in \A$ and
a pair of uniquely determined constituents
$[w_i]_a$ and $[w_i']_a$, for $w_i,w_i'$ otherwise as in (i).
\ere

\medskip
In cases~(ii) and~(iii) the $a\in \A$ is the same in both structures,
or else player~$\PI$ would have won. We think of the $w_i$ and~$w_i'$
in the core parts as the anchors of actual choices made in
corresponding rounds, which may additionally involve an $a_i \in \A$
and subsets $s_i$ and $s_i'$ from the second sorts.
We note that in case~(iii)
these anchor worlds $w_i,w_i'$ are not uniquely determined as representatives of 
their $a_i$-classes; but we shall argue below that matching choices of
these anchors can be made. 
So we describe
the game position after $k$ rounds, for $k \leq r$, in terms of the associated $\wbar =
(w_0,\ldots,w_k)$ and $\wbar' = (w_0',\ldots,w_k')$ and (for some $i$)
additional $a_i$ (as in~(ii)) or $(a_i,s_i)$ and $(a_i,s_i')$ (as in~(iii)). 

We augment the content of the current position further to an invariant
that is to be maintained by player~$\PII$, as follows.

\paragraph*{(II) Abstraction of an invariant.}
With the tuple $\wbar$ of anchor worlds in the core part of $\str{N}$
we associate a tree-like hypergraph structure whose hyperedges represent overlapping 
$a$-classes in the underlying basic modal $S5$-frame. For this we use
a minimal spanning sub-tree containing the worlds listed in $\wbar$. 
We write $\mathrm{tree}(\wbar)$ for this tree structure. 
For its vertices $\mathrm{tree}(\wbar)$ comprises, for every $w_i$ in $\wbar$,
the unique sequence of worlds in which those
$a$-classes intersect that make up the shortest connecting path from the root
$w$ to $w_i$, as well as its end points, $w = w_0$ and
$w_i$. The hyperedges of $\mathrm{tree}(\wbar)$
are formed by the non-trivial subsets of its
vertices that fall within the same $a$-class, labelled
by the corresponding $a$ (any such subset forms a clique w.r.t.\ the 
corresponding $R_a$, for a unique $a \in \A$ due to acyclicity). So
far the $w_i$ are covered by a unique representation of the overlapping $a$-classes 
along these shortest connecting paths to the root. We label each 
vertex $u$ of this tree structure that stems from a world of $\M$ 
by its $\sim^{\ell-d +1}$-type in $\M$, where $d= d(w,u)$ is its  
distance from the root $w = w_0$ (its depth in the tree structure);
and label its hyperedges by the appropriate agents $a \in \A$. 
The tree $\mathrm{tree}(\wbar')$ on the side of $\str{N}'$ is similarly
defined. After $k$ rounds, each of these trees can have at most $1+k\ell$ vertices 
and each individual hyperedge can have at most $k+1$ vertices. 
For worlds $u$ of $\M$ at distance $d = d(w,u)$,
agent $a \in \A$ and $m \in \N$, we now denote as 
$\str{M}\restr [u]_a,\rho_\simm$ the locally full (i.e.\ full) relational
encoding of the local $\simm$-structure $\M \restr [u]_a,\rho_\simm$
at granularity $m = \ell - d +1$, which has colours for
$\simm$-types in
$\M$ (cf. Definition~\ref{localastrgranulardef}); and analogously on
the side of $\str{N}'$.
We argue that player~\PII\ can maintain the following isomorphy
conditions in terms of these tree structures
$\mathrm{tree}(\wbar)$ and $\mathrm{tree}(\wbar')$
through $r$ rounds, and thus force a win.

\medskip
\noindent
$\INV{k}$ after round~$k$:
\bne
\item
the tree structures $\mathrm{tree}(\wbar)$ and $\mathrm{tree}(\wbar')$ 
spanned by $\wbar = (w_0, \ldots,w_k)$ and
$\wbar' = (w_0', \ldots,w_k')$ are isomorphic (as hyperedge- and
vertex-labelled hypergraphs) 
via an isomorphism $\zeta$ that maps $\wbar$ to $\wbar'$:
\[
\zeta \colon \mathrm{tree}(\wbar),\wbar \simeq \mathrm{tree}(\wbar'),\wbar'
\]
\item
  for each pair of worlds $u \in \mathrm{tree}(\wbar)$ and $u' \in \mathrm{tree}(\wbar')$
  that are linked according to $u' = \zeta(u)$ in~(1) at depth $d =
  d(w,u) = d(w',u')$, and for any label $a$:  
\[
  \str{M}\restr [u]_a,\rho_\simm,\sbar
  \equiv_{r-z+1}
  \str{M}'\restr [u']_a,\rho_\simm,\sbar' \quad \mbox{ for } m =\ell-d
  +1, 
\]
where $\sbar$ and $\sbar'$ are tuples of size $z$ that coherently list any 
singleton information states corresponding to the tree vertices incident with that
$a$-hyperedge and any information states $s_i\subset [u]_a = \sigma_a(u)$ and 
$s'_i\subset [u']_a = \sigma_a(u')$ that may have been chosen in the second sort
of the corresponding local component structures during the first $k$
rounds of the game.
\ene

\paragraph*{(III) Maintaining the invariant.}
In the form of $\INV{0}$ these conditions are satisfied at the start of the game: (2)
in this case does not add anything beyond~(1), which in turn is
a consequence of the assumption that 
$\M,w_0  \,\simn\, \M',w'_0$.
We show how to maintain the invariant (conditions~(1) and~(2)) through round $k$,
in which player~$\PI$ may either pebble an element of the first sort,
i.e.\ a world in the core part ($w_k$ or $w'_k$) or in one of the local
component structures ($(w_k)_a$ or $(w'_k)_a$), $a = a_k$, or
an element of the second sort, i.e.\ an information state
($s_k$ in some $[w_k]_a$ or $s_k'$ in some $[w_k']_a$, $a = a_k$,).
We refer to the position before this round by parameters 
$\wbar,\wbar', \mathrm{tree}(\wbar), \dots$ as described above, but at level $k-1$. 
Assuming~$\INV{k-1}$
we need to show how player~$\PII$ can find responses so as to satisfy~$\INV{k}$.
By symmetry we may assume
that $\PI$ opens the $k$-th round by putting the next pebble on the side of $\str{N}$.

\smallskip
Case~1. A move in the second sort: $\PI$ chooses an 
information state $s_k$ so that  $s_k \subset \wp([w_k]_a)$
for some uniquely determined local component $\str{M}\restr [w_k]_a$,
$a = a_k$. 

Case~1.1: if $[w_k]_a$ = $[u]_a$ for some vertex $u$ in 
$\mathrm{tree}(\wbar)$, $u$ at depth $d$ say,
then we look, for $u' = \zeta(u)$, at one round in the game for 
\[
\str{M}\restr [u]_a, \rho_\simm, \sbar\;\equiv_{r-z+1} 
\str{M}'\restr [u']_a, \rho_\simm, \sbar'
\quad \mbox{ for } m = \ell -d +1 
\]
with a move by the first player on $s_k$, which extends $\sbar$ to
$\sbar s_k$; this has an adequate response $s'_k$ for the second
player, which extends $\sbar'$ to $\sbar s'_k$ and 
guarantees $\equiv_{r-z}$ (this is the appropriate level since the tuples
$\sbar$ and $\sbar'$ have been extended by one component). 

Case~1.2: $[w_k]_a$ is ``new'' and the 
appropriate $w'_k$ that satisfies conditions~(1) and~(2) has to be located in a
first step that simulates a move on $(w_k)_a$ (treated as Case~3), after
which we may proceed as in Case~1.1.

\smallskip
Case~2.
Suppose $\PI$ chooses a world in the core part of $\str{N}$ in the
$k$-th round, say $w_k$. The choice of an appropriate match $w_k'$ 
is treated by induction on the distance that the newly chosen 
world $w_k$ has from $\mathrm{tree}(\wbar)$.
In the base case, distance $0$ from $\mathrm{tree}(\wbar)$, 
$w_k$ is a vertex of $\mathrm{tree}(\wbar)$ and nothing needs to
be updated: the response is dictated by the existing isomorphism $\zeta$
according to (1). 
In all other cases, $\mathrm{tree}(\wbar)$ and $\zeta$ need to be
extended to encompass the new $w_k$. The new world~$w_k$ can be 
joined to $\mathrm{tree}(\wbar)$ by a unique shortest path of
overlapping $a$-classes of length greater than~$0$ that connects it to
$\mathrm{tree}(\wbar)$. The new branch in $\mathrm{tree}(\wbar)$ will 
be joined to $\mathrm{tree}(\wbar)$ either through a new $a$-hyperedge
emanating from an existing vertex $u$ of $\mathrm{tree}(\wbar)$ 
(treated in Case~2.1) or through a new vertex $u$ that is part of 
an $a$-class that already contributes an existing hyperedge of  $\mathrm{tree}(\wbar)$
(treated in Case~2.2).

Case~2.1: it is instructive to look at the special case of distance $1$ from
$\mathrm{tree}(\wbar)$ and then argue how to iterate for larger
distance.  
So let $w_k$ be at distance $1$ from 
$\mathrm{tree}(\wbar)$ in the sense that 
$w_k \in R_a[u]$ for some $u$ in
$\mathrm{tree}(\wbar)$ at depth $d$, but $w_k$ not incident with an
$a$-hyperedge in $\mathrm{tree}(\wbar)$ (else in Case~2.2). 
Since $\M, u \sim^{\ell-d+1} \M', u'$ for $u' = \zeta(u)$ and 
\[
\str{M}\restr [u]_a, \rho_\simm,u \equiv_{r} 
\str{M}'\restr [u']_a,\rho_\simm, u' \quad \mbox{ for } m = \ell -d +1,
\footnote{Recall that $\sim^1$ above $\rho_\simm$ suffices for
  Corollary~\ref{localupgradecor}; and $\ell - d+1 = m$.}
\]
a suitable response to the move that pebbles the world $w_k$ (or the singleton
information state $\{w_k\}$) in that game yields a world $w_k' \in E_a[u']$ such that $\M, w_k \sim^{\ell-d} \M', w_k'$
and 
\[
\str{M}\restr [u]_a, \rho_\simm, \{u\}, \{ w_k\} \;\equiv_{r-1} 
\str{M}'\restr [u']_a, \rho_\simm, \{ u'\},\{ w_k'\}
\quad \mbox{ for } m = \ell -d.
\]

These levels of equivalence and granularity are 
appropriate since the depth of $w_k$ and $w_k'$ is $d+1$, and 
since one new vertex contributes to the new hyperedge. 
So we may extend the isomorphism $\zeta$ to map $w_k$ to $w'_k$ 
in keeping with conditions~(1) and~(2); in particular $w_k \sim^{\ell
  - (d+1) +1} w_k'$. If $w_k$ is at greater
distance from its nearest neighbour $u$ in $\mathrm{tree}(\wbar)$ we can
iterate this process of introducing new hyperedges with one new
element at a time, degrading the level of inquisitive bisimulation equivalence by
$1$ in every step that takes us one step further away from the root,
but maintaining equivalences $\equiv_{r-1}$ between the newly added
local $a$-structures 
(recall that $\sim^1$ supports full
bisimilarity of local structures so that this feature does not degrade
with distance, while granularity does). 

Case~2.2: this is the case of
$w_k \in [u]_a$ for some $u$ in $\mathrm{tree}(\wbar)$
that is already incident with an $a$-hyperedge of $\mathrm{tree}(\wbar)$.
By~(2) we have 
\[
\str{M}\restr [u]_a, \rho_\simm, \sbar \;\equiv_{r-z+1} 
\str{M}'\restr [u']_a, \rho_\simm, \sbar' \quad \mbox{ for } m = \ell
-d +1,
\]
where $z$ is the size of the tuples $\sbar$ and $\sbar'$ already
incident with these local component structures, and $d$ is the depth of
$u$ and $u'$. So we can find a response to a move on $\{ w_k\}$ in
this game that yields $w_k' \in [u']_a$ such that 
$\M,w_k \sim^{\ell-d} \M',w'_k$
and  
\[
\str{M}\restr [u]_a, \rho_\simm,\sbar\{w_k\} \;\equiv_{r-z} 
\str{M}'\restr [u']_a, \rho_\simm, \sbar'\{w'_k\}
\quad \mbox{ for } m = \ell -d.
\]

The levels of bisimulation equivalence and granularity are 
appropriate as $m = \ell - (d+1) +1$, the
depth of $w_k$ and $w_k'$ is one greater than that of $u$ and $u'$,
and as the length of the tuples $\sbar$ and
$\sbar'$ has been increased by $1$. So we may extend $\zeta$ by 
matching $w_k$ with $w'_k$ and extending 
$\mathrm{tree}(\wbar)$ and $\mathrm{tree}(\wbar)'$ by these new
vertices and stay consistent with conditions~(1) and~(2).

\smallskip
Case~3.
If $\PI$ chooses some $(w_k)_a$ in the first sort of a local constituent
$\str{M}\restr [w_k]_a$ of $\str{N}$, the choice of an appropriate
match $w_k'$ for $w_k$ is carried out as in Case~2 (with a 
label~$a_k := a$ in the description of the position),
and \PII\ is to respond with a pebble on $(w_k')_a$.
This finishes the proof of Lemma~\ref{inqS5upgradelem} in its essentially
world-oriented formulation. 
\eprf

The adaptations for the non-trivially state-pointed version are marginal.
Recall the illustration of
the upgrading target in Figure~\ref{S5extendedupgradefigure},
equation~(\ref{Gaifmanequiveqn}) and 
Proposition~\ref{summaryGaifmanupgradeprop}. 
By Definition~\ref{explodedviewdef} we use, for the state-pointed case, 
$\expl{\M,s} = \expl{\str{M},s} = \expl{\str{M}_s},s$, based on 
the augmentation $\str{M}_s$ of $\str{M} = \str{M}^\lf(\M)$ 
by the dummy agent $s \in \A_s$ and parameter $s = s_{\alpha_s}$ in the
second sort.  So it suffices to establish
\[
\str{N},s := \expl{\str{M}_s}\restr
N^\ell(s),s  
  \equiv_r^{\ssc(\ell)}
\expl{\str{M}_{s'}}\restr
N^\ell(s'),s' =: \str{N}',s'
\]
for suitable~$\ell$ and~$r$, where $\M_s,s \sim^n \M'_{s'},s'$ for $n \geq \ell +1$, and both
$\M_s$ and $\M'_{s'}$ are $N$-acyclic, now for $N > 2(\ell+2)$, 
sufficiently simple and rich (cf. Proposition~\ref{bisimpartnerprop}).
Here $\M,s \sim^n \M',s'$ implies that for any choice of $w \in s$ we
can find $w' \in s'$ such that $\M,w \sim^n \M',w'$.
We note that $N^\ell(s)\subset N^{\ell+2}(w)$ 
and that $s = s_{\alpha_s}$ as well as 
$N^\ell(s)$ are $\FO$-definable in $\expl{\str{M}_s}\restr
N^{\ell+1}(w),w$ (and similarly on the side of $\M'$). It therefore
suffices to establish $\equiv_{r+t}^{\ssc(\ell+2)}$ 
between $\str{N},w := \expl{\str{M}_s}\restr
N^{\ell+2}(w),w$ and its counterpart $\str{N}',w' := \expl{\str{M}_{s'}}\restr
N^{\ell+2}(w'),w'$ (for a constant $t$ depending on just~$\ell$ that
is large enough to cover the
$\FO$-definitions of $s$ and $N^\ell(s)$ in terms of $w \in s$ in the
exploded views).
Modulo these adaptations of the parameters, this case therefore reduces
to the world-pointed case treated in Lemma~\ref{inqS5upgradelem}.
The initial condition $\INV{0}$ now concerns, on the side of $\str{N}$
say, the anchor $w_0 := w$ in the core, augmented by $a_0 := s \in \A_s$ 
and $s_0 := s_{\alpha_s} = [w]_s$ as induced by the initial pebble on
$s_{\alpha_s}$ in the second sort of the $s$-local component structure
$\str{M}_s\restr [w]_s$ in the exploded view. 

\subsubsection{Summary: back to the main theorem}
Returning to the main theorem, Theorem~\ref{main2}, we recall 
how Corollary~\ref{EFcorr} reduces the $\inqbm$-definability claim
for $\sim$-invariant (and for state properties, inquisitive)
$\FO$-defined world- or state-properties over the respective
classes of relational inquisitive epistemic models to
their $\simn$-invariance, for some sufficiently large finite value for~$n$;
see Observation~\ref{compactobs} to this effect. 
This criterion can
be established through an upgrading of a sufficiently high level of
$\simn$-equivalence to indistinguishability by the given
$\FO$-formula. Through a chain of $\sim$-preserving  model
transformations within the class of models under consideration,
and by passage to the exploded view of these models, this upgrading
eventually boils down to an upgrading of $\simn$
to some level of local Gaifman equivalence, as stated in 
Proposition~\ref{summaryGaifmanupgradeprop}.
Figure~\ref{S5extendedupgradefigure} provides an overview over the several
translation steps in the two tracks for the world-pointed and the
state-pointed version, respectively.
And those upgrading tracks have been completed with
Lemma~\ref{inqS5upgradelem} and its adaptation to the
non-trivially state-pointed case.

So this concludes the proof of the main theorem -- in the
world-pointed and in the state-pointed versions, respectively.
In each case the argument is good in the classical as
well as in the finite model theory reading since all the intermediate
stages in the upgrading argument, within the class of general or
locally full relational encodings, also preserve finiteness. Via 
passage to the exploded views, the upgrading always goes through
locally full relational encodings. The difference between the 
expressive completeness assertions over the class of general versus
locally full relational encodings may therefore almost be overlooked;
but if, for instance $\phi=\phi(\mathtt{s}) \in \FO$ is only assumed to be $\sim$-invariant
and persistent over the class of locally full relational encodings, it can also only be
guaranteed to be $\simn$-invariant for sufficiently large $n$ over
those encodings, and hence equivalently expressible by some
$\phi' \in \inqbm$ over that restricted class of relational encodings;
it is the logical equivalence $\phi' \equiv \phi$ that is subject to
the restriction; and of course $\phi' \equiv \phi$ remains valid over
the class of all relational encodings provided $\phi$ is
$\sim$-invariant throughout, just as $\phi'$ is anyway.

\section{Conclusion}

In our  previous paper~\cite{CiardelliOttoJSL2021} we introduced a
notion of bisimulation for inquisitive modal logic and proved a
counterpart of the van Benthem--Rosen theorem for inquisitive modal
logic over the class of all (finite or arbitrary) relational encodings
of inquisitive modal models.
In this paper, we have extended those results to the especially interesting, but technically challenging, case of inquisitive \emph{epistemic} models, which satisfy $S5$-like constraints. 
We have compared the expressive power of inquisitive modal logic to
that of first-order logic interpreted over relational encodings of
such models, showing that a property of world-pointed inquisitive
epistemic models is expressible in inquisitive epistemic logic just in
case it is expressible in first-order logic and bisimulation
invariant; moreover, the same goes for \emph{inquisitive} properties
of state-pointed models, i.e., properties that are preserved in subsets and always contain the empty state. We have seen that the result holds also if we restrict to the class of finite models; moreover, it holds regardless of whether or not we impose a \emph{local fullness} constraint, which requires the encoding of a model to represent in the second sort the full powerset of each equivalence class in the underlying Kripke frame.

Conceptually, as we discussed in the introduction, the result is
interesting as it implies that inquisitive modal logic is in a precise
sense \emph{expressively complete} as a language to specify
bisimulation-invariant properties of worlds in inquisitive epistemic
models (or bisimulation-invariant inquisitive properties of
information states). Technically,
we deal with first-order definability over non-elementary classes of 
 $2$-sorted relational structures that give some specifically constrained access to
 sets of subsets in their  
 second sort. This takes us into explorations of a border area
 between first-order and monadic second-order expressiveness. As it
 turns out, the model classes in question are still sufficiently malleable up to
 bisimulation to allow for suitably adapted locality arguments in a 
 back-and-forth analysis -- albeit after having translated the
 structures in question into a more manageable regime by means
 of first-order interpretations. In this context an adaptation of 
 well-established classical model-theoretic methods is at
 the core of our treatment; and the specific combination of these
 ingredients can here serve rather intricate non-classical ends. It is
 in this respect that, also technically, the present paper goes 
 substantially beyond~\cite{CiardelliOttoJSL2021}.

Several directions for further work naturally suggest
themselves. First, we may hope to extend this result to other
inquisitive modal logics determined by salient classes of frames; for
instance, to the inquisitive counterpart of the modal logics
\textsf{K4} (determined by the condition: if $v\in\sigma_a(w)$ then
$\Sigma_a(v)\subseteq\Sigma_a(w)$, which implies the transitivity of the underlying Kripke frame) 
and \textsf{S4} (determined by the
above condition plus the requirement $w\in\sigma_a(w)$, which amounts
to the reflexivity of the underlying Kripke frame).
Characterisations of bisimulation-invariant $\FO$ over related 
classes of Kripke structures have been obtained in \cite{DawarOttoAPAL09};  
for discussion of the inquisitive counterparts of some standard frame
conditions, see~\S7.4 of~\cite{IvanoDiss}. 

Moreover, recent work \cite{Ciardelli:22aiml} introduced a generalization of inquisitive modal logic designed to describe models where the set $\Sigma(w)$ of successors of a given world $w$ is not necessarily downward closed. This generalization uses a conditional modality $\Rrightarrow$, where a formula $\phi\Rrightarrow\psi$ is interpreted as claiming that every successor state that supports $\phi$ also supports $\psi$. One interesting question is whether an analogue of our main theorem holds for this generalised inquisitive logic relative to the class of all $S5$-like models, i.e., models where we have $w\in\sigma_a(w)$ for all $w$, and moreover $v\in\sigma_a(w)$ implies $\Sigma_a(w)=\Sigma_a(v)$, but where the set of successors $\Sigma(w)$ is not necessarily downward closed.

\paragraph*{Acknowledgements.}
We are very grateful to two anonymous reviewers for precious
suggestions that led to an improved presentation, and for spotting
a problem in one of the proofs.

\paragraph*{Funding.}
Ivano Ciardelli's research was supported by the German Research
Foundation (DFG, Project number 446711878) and the European Research
Council (ERC, grant agreement number 101116774); Martin Otto's
research was supported by the German Research
Foundation (DFG grant OT~147/6-1: \emph{Constructions and
  Analysis in Hypergraphs of Controlled Acyclicity}.

\end{document}